\documentclass[a4paper,12pt]{amsart}
\usepackage[english]{certus}

\title{Geometric property (T) for box spaces and sofic approximations}

\author{Vadim Alekseev}
\address{Vadim Alekseev, Technische Universit\"{a}t Dresden, Fakultät Mathematik, Institut f\"{u}r Geometrie, 01062, Dresden, Germany}
\email{vadim.alekseev@tu-dresden.de}

\author{Stefan Drigalla}
\address{Stefan Drigalla, Universit{\"a}t Leipzig, Fakultät für Mathematik und Informatik, Mathematisches Institut, PF 10 09 20, 04109 Leipzig, Germany}
\email{sd36lyse@studserv.uni-leipzig.de}


\newcommand{\scp}[1]{\langle #1 \rangle}
\newcommand{\fix}{\mathrm{fix}\,}
\newcommand{\dt}{\,\mathrm{d}t}

\newcommand{\pb}{\partial\beta}

\begin{document}

\onehalfspace

\begin{abstract} 
We prove that every sofic approximation of a property (T) group is approximately isomorphic to one having geometric property (T), and more generally, a box space of graphs which has boundary geometric property (T) is approximately isomorphic to one having geometric property (T). We also prove that a sequence of bounded degree graphs is approximately isomorphic to a disjoint union of expanders  if and only if the Laplacian has spectral gap in the ultraproduct. Finally, we prove a local geometric criterion for geometric property (T) in the spirit of \.{Z}uk's criterion for property (T) for groups.
\end{abstract}

\maketitle

\section{Introduction}
This paper continues the previous work by the first named author with Finn-Sell \cite{AFS} and Biz \cite{alekseev-biz} relating coarse geometry of the graph spaces obtained from sofic approximations of a finitely generated group $\Gamma$ and analytic properties of the group, using the coarse boundary groupoid to connect them. The coarse boundary groupoid can usefully be considered both as topological or a measured groupoid, reflecting the difference between more ``rigid'' coarse geometric properties and ``geometry up to negligible subsets'', usually encountered in the context of sofic approximations.

In \cite{AFS}, it was observed that the coarse properties of the box space of a group that reflect amenability, a-T-menability or property (T) of a residually finite finitely generated group $\Gamma$ can fail for sofic approximations and it was asked which variants of these properties would be the correct replacement in the sofic case. 
This question has been resolved for amenable groups in  \cite{kaiser} and a-T-menable groups in \cite{alekseev-biz}; one of the remaining cases is that of property (T) which we resolve in this article. Moreover, in \cite{alekseev-biz} it was recognized that these results are naturally derived for general sequences of bounded degree graphs and their coarse groupoids, making the statement about sofic approximations their direct corollaries. It turned out that measurable amenability or a-T-menability of the coarse boundary groupoid corresponds to the space of graphs having the expected coarse property (property A resp. asymptotic coarse embeddability into a Hilbert space) after possibly changing the graphs on subsets which have measure zero in the ultralimit (thus passing to what we will call an approximately isomorphic sequence).

In this article we therefore take the box space of a bounded degree graph sequence as our primary object of study. In the work of Willett and Yu \cite{WYgeomT} geometric property (T) for box spaces has been identified as the correct counterpart to Kazhdan's property (T) for a group $\Gamma$ in the residually finite case; it also has an interpretation through property (T) of the coarse groupoid by \cite{dA-W}.  On the other hand, by the work of Kun \cite{kun} every sofic approximation of a property (T) group is approximately isomorphic to a sequence consisting of (disjoint unions of) expander graphs; being an expander\footnote{In this paper, we do not require graphs constituting an expander sequence to be connected.} sequence is known to be strictly weaker than having geometric property (T). Therefore it was natural to expect that the coarse boundary groupoid having property (T) on a subset of measure $1$ with respect to an ultralimit of counting measures should imply that the graph sequence is approximately isomorphic to one which actually has geometric property (T), and our first main result confirms this:

\begin{ThmA}[{Theorem \ref{theorem: almost T implies appr eq to geometric T}}]
Let $X$ be a box space of graphs which has almost boundary geometric property (T). Then there is a box space $X'$ that is approximately isomorphic to $X$ and that has geometric property (T).
\end{ThmA}

Similarly to \cite{alekseev-biz}, it implies the corresponding statement for sofic approximations:
\begin{ThmA}[{Theorem \ref{thm: almost T sofic approx}}]
Let $\Gamma$ be a finitely generated sofic group. The following statements are equivalent:
\begin{enumerate}
    \item $\Gamma$ has property (T),
    \item every sofic approximation of $\Gamma$ has almost boundary geometric property (T),
    \item there is a sofic approximation of $\Gamma$ with geometric property (T).
\end{enumerate}
Moreover, every sofic approximation of $\Gamma$ is approximately isomorphic to one having geometric property (T).
\end{ThmA}

Since the coarse groupoid of $X$ decomposes into the ``local'' part (the restriction to $X$) and the ``boundary part'' (the restriction to the Stone-\v{C}ech boundary of $X$), it turns out that geometric property (T) also decomposes into having spectral gap for the Laplacian both on $\ell^2(X)$ (that is, being an expander) and on the Stone-\v{C}ech boundary of $X$; the latter property is then called \emph{boundary property (T)}. Therefore an important part of the proof consists in deducing that $X$ is approximately isomorphic to an expander sequence. In fact, we were able to prove a complete characterization of this property through the spectral gap in the ultraproduct.

\begin{ThmA}[{Corollary \ref{corollary:almost T is appr expander}}]
A box space of graphs $X$ has spectral gap in the ultraproduct along every ultrafilter $\mk u\in \partial\beta \mb N$ if and only if $X$ is approximately isomorphic to an expander. 
\end{ThmA}

In the proof of this theorem we use an adaptation of the method of Kun introduced in \cite{kun} in context of sofic approximations, adapting it both to sequences of bounded degree graphs and to the weaker assumption of just having spectral gap in the ultraproduct and not the full force of property (T) which gave the possibility to apply spectral gap to arbitrarily small subsets in the graph sequence in \cite{kun}.

Similarly to \cite{alekseev-biz}, it is natural to put our considerations into context of discrete p.m.p. measured groupoids. To this end, we introduce a notion of property (T) for discrete p.m.p. measured groupoids which works well in the context of the coarse boundary groupoid equipped with an ultralimit of counting measures and prove that for standard Borel groupoids it coincides with the definition from \cite{A-D}. On the way, we deduce some results about property (T) for measured groupoids in parallel to the results for topological groupoids from \cite{dA-W}. In particular, a compactly generated étale topological groupoid with (topological) property (T) with an invariant measure on the base space has measured property (T), and the latter implies spectral gap in the ultraproduct.

Finally, we prove a ``local'' geometric criterion for geometric property (T) in the spirit of \.{Z}uk's criterion for property (T) for groups:
\begin{ThmA}[{Theorem \ref{thm:geomZuk}}]
    Let $X=\bigsqcup_i X_i$ be a box space of graphs such that each linking graph $L_x$ is connected. Suppose that there are $Y_i\subseteq X_i$ with $|Y_i|/|X_i|\to 1$ and $\lambda > \frac{1}{2}$ such that for every $x\in Y_i$ we have 
    \[
    \lambda_1\left(L_x\right)\geq \lambda.
    \]
    Then $X$ has almost boundary (T) with respect to $\mathrm{core}(Y)$. If $Y_i = X_i$ for all $i\in\mathbb N$, then $X$ has geometric property (T).
\end{ThmA}

\textbf{Acknowledgements.} The authors would like to thank Rufus Willett for comments and suggestions that helped to improve the text, and Andreas Thom for stimulating discussions. The authors acknowledge support by the DFG (SPP 2026 ``Geometry at infinity''). Some of the results were obtained as part of the second namend author's Ph.D. thesis.

\section{Preliminaries}

We are interested in the coarse properties of sequences of bounded degree graphs. The main objects of study for us are the following:

\begin{Def}
A \emph{box space of graphs} is the disjoint union $X=\bigsqcup_i X_i$, where $(X_i)_{i\in\mb N}$ is a sequence of finite graphs of uniformly bounded degree satisfying $|X_i|\to \infty$. We equip $X$ with a metric that restricts to the graph metric on each connected component of $X_i$ and is $\infty$ else. We do \textbf{not} require the $X_i$ to be connected for reasons that will become apparent later. 
\end{Def}

\begin{Rem}
To avoid using infinite-valued metrics, one could alternatively choose a metric on each $X_i$ which restricts to the graph metric on each connected component and such that
\[
d(X_i^{(k)}, X_j^{(k')})\to \infty
\]
as $i+j\to \infty$, where $X_i^{(k)}$, $X_j^{(k')}$ are connected components of $X_i$ resp. $X_j$. See also Remark 2.5 and Remark 3.6 in \cite{WYgeomT}.
\end{Rem}

Although it would be possible to put some of the results in the setting of abstract coarse spaces, we believe it to be more natural to stick with the metric setting entirely. It is worth noting that under reasonable assumptions, e.g. those made in \cite{WYgeomT}, a coarse space is coarsely equivalent to a metric space (possibly with infinite distances), see \cite[Section 2.4]{Roe}.

There are two common approaches of investigating the coarse properties of box spaces: translation algebras and coarse groupoids. As we will want to make use of both, we will gather here the necessary background:

\begin{Def}
         We define the \emph{algebraic Roe algebra} $\mb C[X]$ to be the set of all $T\in \mc B(\ell^2(X))$, such that there is $R>0$ and $\ip{T\delta_x,\delta_y}=0$ for all $x,y\in X$ with $d(x,y)\geq R$. This naturally forms a $\ast$-subalgebra. Alternatively, we can view the elements of $\mb C[X]$ as matrices $(T_{xy})$ indexed by $X\times X$ and values in $\mb C$ through $T_{xy}=\ip{T\delta_x,\delta_y}$. By a representation of any unital $\ast$-algebra we mean a unital $\ast$-homomorphism $\pi:\mb C[X]\to\mc B(\mc H)$ for some Hilbert space $\mc H$. We are particularly interested in the closure of $\mb C[X]$ with respect to the norm
\[
\norm{T}\coloneqq\sup\left\{\norm{\pi(T)}_{\mc B(\mc H)}\mid \pi:\mb C[X]\to\mc B(\mc H) \text{ a representation}\right\} 
\]
denoted $C^*_{u,\max}(X)$, the \emph{maximal uniform Roe algebra} of $X$.
\end{Def}

\begin{Def}
Given a box space of graphs $X=\bigsqcup_i X_i$, or any (discrete) topological space, we denote by $\beta X$ its \emph{Stone-\v{C}ech compactification} and we call $\partial \beta X=\beta X\setminus X$ its \emph{Stone-\v{C}ech boundary}. Fixing any non-principal ultrafilter $\mk u\in\partial\beta\mb N$ we obtain a measure on $\beta X$ corresponding to the state
\[
\mu_{\mk u}(f)=\lim_{i\to\mk u}\frac{1}{|X_i|}\sum_{x\in X_i}f(x),\quad f\in C(\beta X).
\]
\end{Def}

\begin{Rem}
Since $\mu_{\mk u}(X_i)=0$ for all $i$, we also have $\mu_{\mk u}(X)=0$ and thus $\mu_{\mk u}(\partial\beta X)=\mu_{\mk u}(\beta X)=1$. We will also need the fact that $\mu_{\mk u}$ is a regular Borel measure.
\end{Rem}

\begin{Def}
Given a non-principal ultrafilter $\mk u\in\partial\beta\mb N$, we will say that a condition holds along $\mk u$ if there is $I\in\mk u$ such that the condition holds for all $i\in I$.
\end{Def}

\begin{Def}
    A \emph{groupoid} is a set $G$ equipped with the source and range maps $s,r:G\to G^{(0)}\subseteq G$ and a partially defined composition $(g,h)\in G^{(2)}\mapsto gh\in G$, where $G^{(2)}=\{(g,h)\in G\times G:\, s(g)=r(h)\}$, such that for every $x\in G^{(0)}$ and $g,h,k\in G$ with $s(g)=r(h)$ and $s(h)=r(k)$ we have:
    \begin{enumerate}
        \item $s(gh)=s(h)$ and $r(gh)=r(g)$,
        \item $(gh)k=g(hk)$,
        \item $x=s(x)=r(x)$,
        \item $g=gs(g)=r(g)g$,
        \item there is $g^{-1}\in G$ such that $g^{-1}g=s(g)$ and $gg^{-1}=r(g)$.
    \end{enumerate}
    The subset $G^{(0)}$ is called the \emph{base space} or \emph{space of objects} and the elements of $G^{(2)}$ are called \emph{composable pairs}. If $G$ carries a topology so that source, range, composition and inverse are continuous, then $G$ is a \emph{topological groupoid}. It is further called \emph{\'etale} if $s$ and $r$ are local homeomorphisms. If $G$ carries a Borel structure making source, range, composition and inverse measurable maps, then $G$ is a \emph{Borel groupoid}.
\end{Def}

\begin{Def}
For $R>0$ we define $E_R=\left\{(x,y)\in X\times X\colon\ d(x,y)\leq R\right\}$ and 
\[
G(X)=\bigcup_{R>0}\overline{E_R}\subseteq\beta X\times \beta X
\]
the \emph{coarse groupoid of X}, where the closure is taken in $\beta X\times \beta X$. With the weak topology coming from the union and groupoid operation coming from the pair groupoid $\beta X\times \beta X$ this becomes a Hausdorff \'etale topological groupoid with base space $\beta X$, compare \cite[Chapter 10]{Roe}. 
\end{Def}

\begin{Def}
A subset $Z\subseteq G^{(0)}$ is called \emph{invariant} or \emph{saturated}, if for every $g\in G$ we have $s(g)\in Z$ if and only if $r(g)\in Z$. The set $G|_Z$ of all $g$ for which this holds true then becomes a subgroupoid of $G$, which is called the \emph{reduction} of $G$ to $Z$. We will denote the \emph{boundary groupoid} $G(X)|_{\partial \beta X}$ by $\partial G(X)$.
\end{Def}

In what follows all topological groupoids we are going to consider will be Hausdorff and \'etale.

\begin{Def}
Let $G$ be a topological groupoid. Consider the convolution and involution
\[
(f_1\ast f_2)(g)=\sum_{g=hk} f_1(h)f_2(k)\quad \text{and}\quad f^*(g)=\overline{f(g^{-1})}
\]
for $f,f_1,f_2\in C_c(G)$. This turns $C_c(G)$ into $\ast$-algebra, called the \emph{groupoid algebra} of $G$. Its completion with respect to the norm
\[
\norm{T}\coloneqq\sup\left\{\norm{\pi(T)}_{\mc B(\mc H)}\mid \pi: C_c(G)\to\mc B(\mc H) \text{ a representation}\right\} 
\]
is called the maximal groupoid C*-algebra of $G$, denoted by $C^*_{\max}(G)$.
\end{Def}

\begin{Rem}
There is an isomorphism of $\ast$-algebras between the algebraic Roe algebra $\mb C[X]$ and the groupoid algebra $C_c(G(X))$, that is $\mb C[X] \cong C_c(G(X))$, which naturally extends to an isometric isomorphism
\[
C^*_{u,\max}(X) \cong C^*_{\max}G(X).
\]
\end{Rem}

In Section \ref{section: measured T} we will also consider discrete probability measure preserving (p.m.p.) measured groupoids.  We recall only the basic definitions and notation here referring to \cite{sauer, amenable-groupoids} for the general theory of discrete measured groupoids. 

\begin{Def}
A discrete p.m.p. measured groupoid is a Borel groupoid $G$ together with a probability measure $\mu$ on the base space $G^{(0)}$ such that the measure $\mu_{G}$ on $G$ defined on a Borel subset $A\subseteq G$ as
\beqn
\mu_{G}(A):= \int_{G^{(0)}} |s^{-1}(x)\cap A|\, d\mu(x)
\eeqn
coincides with the measure
\beqn
\mu'_{G}(A):= \int_{G^{(0)}} |r^{-1}(x)\cap A|\, d\mu(x)
\eeqn
(such a measure $\mu$ on $G^{(0)}$ is called $G$-invariant).

A measurable subset $K\subseteq G$ is called \emph{bounded} if there exists an $M>0$ such that for almost every $x\in G^{(0)}$ both $|s^{-1}(x) \cap K| < M$ and $|s^{-1}(x) \cap K^{-1}| < M$.  Let $\mb CG$ be the groupoid $\ast$-algebra of $G$, defined as
\beqn
\mb CG:=\{f\in L^\infty(G,\mu_{G})\,|\, \supp(f)\text{ is a bounded subset of }G\}.
\eeqn
with the convolution and involution defined as
\[
(f_1\ast f_2)(g)=\sum_{g=hk} f_1(h)f_2(k)\quad \text{and}\quad f^*(g)=\overline{f(g^{-1})}
\]
The groupoid $\ast$-algebra of $G$ is equipped with a natural linear map
\[
\Psi\colon \mb CG \to L^\infty(G^{(0)}),\quad \Psi(f)(x)= \sum_{g\in r^{-1}(x)} f(g)
\]
\end{Def}

If $G$ is an Hausdorff étale topological groupoid equipped with an invariant probability measure $\mu$ on $G^{(0)}$, then it is straightforward to check that $(G,\mu)$ is a discrete p.m.p. measured groupoid and we get an inclusion of groupoid algebras $C_c(G)\subseteq \mb CG$. In particular, the coarse groupoid $G(X)$ of a box space $X=\bigsqcup_i X_i$ equipped with the measure $\mu_{\mk u}$ defined above is a discrete p.m.p. measured groupoid.

We recall the connection of the coarse groupoid to ultralimits and ultraproducts, as explained in \cite{AFS}. The natural projection map $p\colon X\to \mb N$ extends to $\beta p\colon \beta X\to \beta \mb N$ and we have $(\beta p)(\partial \beta X)=\partial \beta \mb N$. Given a non-principal ultrafilter $\mk u\in \partial \beta \mb N$, the fiber $(\beta p)^{-1}(\mk u)$ is $\partial G(X)$-invariant, and $\mu_{\mk u}$ is supported on it; as a measure space $(\partial G(X)_{\mk u},\mu_{\mk u})$ can then be identified with the measure space ultraproduct $\prod_{i\to \mk u}(X_i,\mu_i)$. We refer to \cite{AFS} for more details and connections of $\partial G(X)|_{(\beta p)^{-1}(\mk u)}$ to graph ultralimits (along $\mk u$) of the graph sequence $(X_i)$.

Altogether, we thus get an isomorphism of von Neumann algebras
\[
L^\infty(\partial\beta X,\mu_{\mk u}) \cong \prod_{i\to \mk u} L^\infty(X_i,\mu_i).
\]

As pointed out in \cite[Section I.3]{connes76}, the ultraproduct von Neumann algebra $L^\infty(\partial\beta X,\mu_{\mk u})$ \emph{does not act} on the ultraproduct Hilbert space $\prod_{i\to \mk u} L^2(X_i,\mu_i)$; instead, the correct Hilbert space on which it acts (the GNS spaces of the state $\mu_{\mk u}$) can be identifeid with a proper subspace of $\prod_{i\to \mk u} L^2(X_i,\mu_i)$ consisting of uniformly integrable sequences \cite[Proposition I.3.1]{connes76}. We refer to \cite{connes76} for more details and, following Connes, denote this ``restricted ultraproduct'' Hilbert space by $\leftsub{{i\to\mk u}}{\prod} L^2(X_i,\mu_i)$. In particular, we get a natural $\ast$-representation of the measured groupoid algebra $\mb C(\partial G(X),\mu_{\mk u})$ on $\leftsub{{i\to\mk u}}{\prod} L^2(X_i,\mu_i)$.

Our focus regarding anything related to property (T) will lie on the {Laplace operator} which is defined as follows:
\begin{Def}
The \emph{Laplace operator} $\Delta\in\mb C[X]$ is given by
\[
(\Delta)_{xy}=\begin{cases}
    \deg (x) &\text{if }x=y, \\
    -1 &\text{if }x\sim y,\\
    0 &\text{else }
\end{cases}
\]
for $x,y\in X$, or equivalently 
\[
(\Delta f)(x)=\sum_{y\sim x}(f(x)-f(y))
\]
for $f\in\ell^2(X)$.
\end{Def}

If $A$ is some $C^*$-algebra and $a\in A$ we denote the spectrum of $a$ by $\sigma_A(a)$. We will also write $\sigma_{\max}(a)$ if $A=C^*_{\max}(G(X))\cong C^*_{u,\max}(X)$ and $\sigma_{\max,Z}(a)$ if $A=C^*_{\max}(\partial G(X)|_Z)$ for $Z\subseteq\partial\beta X$.

\begin{Def}
Let $X=\bigsqcup_i X_i$ be a box space of graphs. We say that $X$ is an \emph{expander} if there is a $c>0$ such that $\sigma_{\mc B(\ell^2(X))}(\Delta)\subseteq\{0\}\cup [c,\infty)$. Equivalently we can require that there is a $c>0$ such that each $X_i$ a union of graphs with Cheeger constant at least $c$. As with general box spaces, we do \textbf{not} require the $X_i$ to be connected, as is often done in the literature.

We now define $X$ to have \emph{geometric property (T)} if there is a $c>0$ such that $\sigma_{\max}(\Delta)\subseteq\{0\}\cup [c,\infty)$. If $\sigma_{\max,\partial \beta X}(\Delta)\subseteq\{0\}\cup [c,\infty)$ then $X$ is said to have \emph{boundary geometric property (T)}.
\end{Def}

\begin{Rem}
    More generally, it holds that a coarse space $X$ has geometric property (T) if $\Delta^E$ has spectral gap in $C^*_{u,\max}(X)$ for some (and equivalently, any) controlled set $E\subseteq X\times X$, see \cite[Theorem 5.8]{WYgeomT}. Here $\Delta^E$ is defined by
    \[
    (\Delta^E)_{xy}=\begin{cases}
    |\{z\in X\colon\;(x,z)\in(E\cup E^{-1})\setminus\diag(E)\}| &\text{if }x=y, \\
    -1 &\text{if }(x,y)\in(E\cup E^{-1})\setminus\diag(E),\\
    0 &\text{else. }
\end{cases}
    \]
    The usual graph Laplacian is obtained for 
    \[
    E=\{(x,y)\in X\times X\colon\;d(x,y)\leq 1\}
    \]
    in case the coarse structure is induced by the metric. In section \ref{zuksection} we will introduce yet another type of Laplacian in order to show a sufficient condition for geometric property (T).
\end{Rem}

So far all known examples of geometric property (T) spaces come from approximations of residually finite groups by a sequence of graphs coming from finite quotients, see \cite[Chapter 7]{WYgeomT}. However, as done in \cite{AFS}, in this context it is natural to consider a broader class of groups:

\begin{Def}\label{defsofic}
    A discrete group $\Gamma$ is called \emph{sofic} if for every finite $F\subseteq\Gamma$ and every $\varepsilon>0$, there is a set $X$ together with a map $\sigma\colon F\to \Sym(X)$ and a subset $Y\subseteq X$ such that:
    \begin{enumerate}
        \item $|Y|\geq (1-\epsilon)|X|$,
        \item $\sigma(g)\sigma(h)(y)=\sigma(gh)(y),\quad g,h\in F,y\in Y$,
        \item $\sigma(g)(y)\neq y,\quad g\in F\setminus\{e\},y\in Y$.
    \end{enumerate}
    Choosing a sequence $(F_i)$ with $F_i\subseteq F_{i+1}$ and $\Gamma=\bigcup_iF_i$ and a sequence $\varepsilon_i\to 0$ the corresponding sequence $(X_i)$ is called a \emph{sofic approximation} of $\Gamma$. 
    If $\Gamma$ is generated by a finite, symmetric set $S$ we obtain a graph structure on each $X_i$ by letting $x \sim y$ whenever there is $s\in S$ with $\sigma_i(s)(x)=y$, and we will always treat $(X_i)$ as a sequence of graphs. Our definition of soficity is then equivalent to saying that these (labelled) graphs converge to the Cayley graph of $\Gamma$ in the sense of Benjamini-Schramm, see \cite[Theorem 5.1]{Pestov}.
\end{Def}

We remark that soficity generalises both being residually finite (and, more generally, being LEF -- locally embeddable into finite groups). We also note that there is, at time of writing, no group that is known to be non-sofic.

\begin{Ex}
    Consider a finitely generated and residually finite group $\Gamma$ and set $X_i=\Gamma/N_i$ be the corresponding sequence of finite quotient groups as in \cite[Definition 7.2]{WYgeomT}. Then it is in particular a sofic approximation and \cite[Theorem 7.3]{WYgeomT} states that $\Gamma$ has property (T) if and only if $X=\bigsqcup_iX_i$ has geometric property (T). The same is true for $\Gamma$ an LEF group, see \cite[Theorem 4.27]{AFS}.
\end{Ex}

\begin{Ex}\label{almostconnectedapproximation}
Take $\Gamma$ with property (T), e.g. $SL(n,\mb Z)$ for $n\geq 3$, and $X$ as above. We obtain a connected graph $X_i'$ by adding an edge between (two arbitrary vertices of) two copies of $X_i$. The sequence $(X_i')$ is still a sofic approximation of $\Gamma$ and the resulting box space $X'$ then has almost boundary (T). Now, since $X'$ fails to be an expander, it does not have geometric property (T). However, $X'$ is by construction approximately isomorphic to a (non-connected) expander, that in fact has geometric property (T). The results in the following sections will show that the general case behaves similarly.
\end{Ex}

\begin{Def}
    Let $X=\bigsqcup_i X_i$ be a box space of graphs with subsets $Y_i\subseteq X_i$ satisfying $|Y_i|/|X_i|\to 1$. Denote $Y=\bigcup_i Y_i$ and define its \emph{core} $Z$ by
    \[
    Z=\mathrm{core}(Y)=(\mathrm{sat}(\pb Y^c))^c\subseteq\pb X,
    \]
    where $\mathrm{sat}(A)=\{x\in \beta X\,:\,(a,x)\in G(X)\text{ for some }a\in A\}$ is the \emph{saturation} of a subset $A\subseteq \beta X$. The box space $X$ is then said to have \emph{almost boundary  geometric property (T)} (or \emph{almost boundary (T)} for short), if $\Delta$ has spectral gap in $C^*_{\max}(G(X)|_Z)$, that is we find $c>0$ such that 
    \[
    \sigma_{C^*_{\max}(G(X)|_Z)}(\Delta)\subseteq\{0\}\cup [c,\infty).
    \]
\end{Def}

The core is invariant by definition and consist precisely of all the points with infinite distance from the complement of $\pb Y$. We formulate a lemma to gather some more useful facts.

\begin{Lemma}\label{coreproperties}
    The core $Z=\mathrm{core}(Y)$ of a box space $X$ is closed, has full measure and satisfies
    \[
    Z=\bigcap_{R>0}\pb (B_R(Y^c))^c.
    \]
\end{Lemma}

\begin{proof}
    Using the notation from the definition we will show that $\mathrm{sat}(\pb Y^c)$ is an open set of measure zero. We clearly have
    \[
    \mathrm{sat}(\pb Y^c)=s(r^{-1}(\pb Y^c))=\bigcup_{R>0}s(\overline{E_R}\cap(\pb Y^c\times \pb X)).
    \]
    Taking into account that $s$ is an open map, $\overline{E_R}$ is clopen and $\pb Y^c$ is open subset of the boundary we see that $\mathrm{sat}(\pb Y^c)$ is open. Further rewriting $$s(\overline{E_R}\cap(\pb Y^c\times \pb X))=\pb (B_R(Y^c))$$ it remains to show that these sets have measure zero for every $R>0$. Since the $X_i$ have uniformly bounded degree, say $d$, we can find $C>0$, only depending on $R$ and $d$, with $|B_R(Y_i^c)|\leq C|Y_i^c|$. Thus we obtain
    \[
    \frac{|B_R(Y_i^c)|}{|X_i|}\leq C\frac{|Y_i^c|}{|X_i|}= C\left(1-\frac{|Y_i|}{|X_i|}\right)\to 0,
    \]
    which is to say that $\overline{B_R(Y^c)}$ has measure zero.
\end{proof}

With the description of $Z$ from the lemma one easily sees that in the case of a sofic approximation the definition of the core is equivalent to the one of the sofic core in \cite[Definition 3.6]{AFS} and furthermore that:

\begin{Cor}\label{soficalmostboundaryT}
    Let $\Gamma$ be a finitely generated group and $(X_i)$ be a sofic approximation. Let $Z=\mathrm{core}(Y)$ be a core of $X=\bigsqcup_i X_i$. Then $\Gamma$ has property (T) if and only if $X$ has almost boundary (T) with respect to $Z$.
\end{Cor}

The next aspect we want to transfer into the setting of box spaces, is that sofic approximations behave well under small changes. We employ the following definition, that is also used in \cite{Winkel}.

\begin{Def}
    Let $X=\bigsqcup_i X_i$ and $X'=\bigsqcup_i X_i'$ be box spaces of graphs. If there are subgraphs $Y_i\subseteq X_i$ and $Y_i'\subseteq X_i'$ such that $Y_i$ and $Y_i'$ are isomorphic as graphs and
    \[
    \lim_{i\to\infty}\frac{|Y_i|}{|X_i|}=\lim_{i\to\infty}\frac{|Y_i'|}{|X_i'|}=\lim_{i\to\infty}\frac{|E(Y_i)|}{|E(X_i)|}=\lim_{i\to\infty}\frac{|E(Y_i')|}{|E(X_i')|}=1,
    \]
    then $X$ and $X'$ are said to be \emph{approximately isomorphic}.
\end{Def}

\begin{Rem}\label{apprisomVScoarseequivalence}
    Obviously, being approximately isomorphic is an equivalence relation that preserves sofic approximations. It does however \textbf{not} imply being coarsely equivalent, as Example \ref{almostconnectedapproximation} tells us. In particular, expansion and geometric property (T) are generally not preserved when passing to an approximately isomorphic sequence.
\end{Rem}

Our next aim is to show that almost boundary (T) is indeed an invariant of approximate isomorphisms. The key observation for this will be formulated in a lemma.

\begin{Lemma}
    Let $X=\bigsqcup_i X_i$ and $X'=\bigsqcup_i X_i'$ be box spaces of graphs that are approximately isomorphic. Then there are isomorphic subgraphs $V_i\subseteq X_i$ and $V_i'\subseteq X'_i$ such that $|V_i|/|X_i|\to 1$, $|V'_i|/|X'_i|\to 1$ and for $W=\mathrm{core}(V)$, $W'=\mathrm{core}(V')$ we have $\partial G(X)|_W \cong \partial G(X')|_{W'}$.
\end{Lemma}

\begin{proof}
Consider sequences of subgraphs $U_i\subset X_i$ and $U'_i\subset X_i'$ from the definition of the approximate isomorphism. for the sake of readability we will write $U_i$ for either of these graphs and consider them as subsets of both $X_i$ and $X_i'$. Define $V_i\subset U_i$ to be the maximal subgraph that has no vertex adjacent to $X_i\setminus U_i$ or $X_i'\setminus U_i$. Since
\[
|V_i|\geq |U_i|-d(|X_i\setminus U_i|+|X_i'\setminus U_i|),
\]
where $d$ is a bound for the degree of $X$ and $X'$, we have
\[
\lim_i \frac{|V_i|}{|X_i|}=\lim_i \frac{|V_i|}{|X_i'|}=1
\]
and the same is true for the edge sets $E(V_i)$. Denote $U=\bigcup_i U_i$ and $V=\bigcup_i V_i$ and define $W=\mathrm{core}(V)$. This definition is unambiguous: First, since $\overline{U}$ is homeomorphic to $\beta U$ we see that $\pb V$ does not depend on the ambient space. Second, by Lemma \ref{coreproperties} we have $W=\bigcap_{R>0}\pb (B_R(V^c))^c$ and as $B_R(V^c)^c\subseteq U$ for all $R>0$ considered as subsets of either $X$ or $X'$ it follows that $W\subseteq\pb U$ is also independent of the ambient space. 

Now write 
\begin{align*}
    E_R^Y &= E_R\cap(Y\times Y)=\{(x,y)\in Y\times Y\,:\,\rho_X(x,y)\leq R\} \subseteq X\times X,\\
    F_R^Y &= F_R\cap(Y\times Y)=\{(x,y)\in Y\times Y\,:\,\rho_{X'}(x,y)\leq R\} \subseteq X'\times X'
\end{align*}
for any subset $Y$ of $X$ or $X'$ respectively. If we further write $V^R=B_R(V^c)^c$ we see that
\[
V^R=\{x\in V\,:\,\rho_X(x,X\setminus V)\geq R\}=\{x\in V:\, \rho_{X'}(x,X'\setminus V)\geq R\},
\]
since any path out of $V_i$ does so through $U_i$. Similarly we see that for $\frac{R}{4}<S<\frac{R}{2}$ we must have $E^{V^R}_S=F^{V^R}_S$. Since we also have
\[
W=\bigcap_{R>0}\pb (B_R(V^c))^c= \bigcap_{R>0}\pb V^R
\]
it follows that $\pb E^W_S\subseteq\pb E^{\pb V^{3S}}_S=\pb F^{\partial V^{3S}}_S\supseteq\pb F^{W}_S$, from which 
\[
\pb E^W_S=\pb E^{\pb V^{3S}}_S\cap (W\times W)=\pb F^{W}_S
\]
follows. This yields
\[
\partial G(X)|_W=\bigcup_{S>0}\pb E^W_{S}= \partial G(X')|_W
\]
and this finishes the proof, because the coarse groupoid carries the weak topology coming from the union of these sets.
\end{proof}

With this at hand we arrive at the desired invariance.

\begin{Prop}\label{almTinv}
    Let $X=\bigsqcup_i X_i$ and $X'=\bigsqcup_i X_i'$ box spaces of graphs that are approximately isomorphic. Then if $X'$ has almost boundary (T), so does $X$.
\end{Prop}

\begin{proof}
    Since $X'$ has almost boundary (T) we find $Y'_i\subseteq X'_i$, $Y'=\bigcup_i Y'_i$ with 
    \[
    |Y'_i|/|X'_i|\to 1\quad \text{and}\quad Z'=\mathrm{core}(Y')
    \]
    such that $\Delta \in C^*_{\max}(\partial G(X')|_{Z'})$ has spectral gap. Take now $V=\bigcup_i V_i$ as in the preceding lemma, then denote $Y_i=Y'_i\cap V_i$ and $Y=\bigcup_i Y_i$. It follows that
    \begin{align*}
        \overline{Y} &= \overline{\bigcup_i (Y'_i\cap V_i)} = \overline{\bigcup_i Y'_i}\cap \overline{\bigcup_i V_i}=\overline{Y'}\cap \overline{V},
    \end{align*}
    where we used that $\overline{A\cap B}=\overline{A}\cap\overline{B}$ for subsets $A,B\subseteq X$ (this follows from the identification $C(\beta X)\cong \ell^\infty X$ by which $\chi_{\overline{A}}\chi_{\overline{B}}$ and $\chi_{\overline{A\cap B}}$ both agree with the projection $p_A p_B=p_{A\cap B}$). Also we clearly have $|Y_i|/|X_i|\to 1$. Thus for $Z=\mathrm{core}(Y)$ we see that restriction from $\partial G(X')|_{Z'}$ to $\partial G(X')|_{Z}$, that by the the preceding lemma can be identified with $\partial G(X)|_{Z}$, induces a map $C^*_{\max}(\partial G(X')|_{Z'})\to C^*_{\max}(\partial G(X)|_{Z})$. By extend $\Delta\in C^*_{\max}(\partial G(X)|_Z)$ will inherit the spectral gap, which proves almost boundary (T) for $X$.
\end{proof}

\section{Almost boundary (T) and expansion}\label{section:almonst T and expansion}
Now that we have established almost boundary (T), we turn to the relation between almost boundary (T), expansion and geometric property (T). The main idea is that geometric property (T) should split into expansion plus a part ``at infinity". Before giving the first result in this direction, we state a short lemma.

\begin{Prop} \label{exp+bdr}
    Let $X = \bigsqcup_i X_i$ be a box space of graphs. Then the following are equivalent:
    \begin{enumerate}
    \item[(i)] $X$ has geometric property (T),
    \item[(ii)] $X$ is an expander and has boundary property (T).
    \end{enumerate}
    Moreover, the size of the spectral gap of the Laplacian is preserved.
\end{Prop}

For the proof we will need the following lemma:
\begin{Lemma}
    Let $X$ be a countable set and $G\subseteq X\times X$ be a subgroupoid of the pair groupoid. Then $C^*_{\max}(G)\subseteq \mathcal K(\ell^2X)$ with equality if $G=X\times X$.
\end{Lemma}

\begin{proof}
    Since a compactly supported function on $X\times X$ must indeed have finite support, an element of $C_c(G)$ is just a finite-rank operator on $\ell^2X$ and thus a compact operator. Consider $C^*_{\max}(G)\to \mathcal K(\ell^2X)$ the extension of this inclusion. Now, each $a\in C_c(G)$ lies in some finite-dimensional matrix algebra which has the norm of $\mathcal K(\ell^2X)$ as its unique norm, so we have $\norm{a}_{C^*_{\max}(G)}=\norm{a}_{\mathcal K(\ell^2X)}$. Thus the above map is isometric and we can identify $C^*_{\max}(G)$ with a C*-subalgebra of $ \mathcal K(\ell^2X)$. Since $C_c(X\times X)$ is all the finite-rank operators, the second statement also follows.
\end{proof}

\begin{proof}[Proof of Proposition \ref{exp+bdr}]
    As already stated before, geometric property (T) implies both expansion and boundary property (T), so we are left to show $(ii)\Rightarrow (i)$. Consider the the short exact sequence coming from the restriction to the boundary $\pb X$:
    \begin{equation*}
    0 \to C^*_{\max}(G(X)|_X) \to C^*_{\max}(G(X)) \overset{q}{\to} C^*_{\max}(\partial G(X)) \to 0.
    \end{equation*}
    First observe that we have $C^*_{\max}(G(X)|_X)\subseteq \mathcal K(\ell^2 X)$: since $G(X)|_X\subseteq X\times X$. Let $\Delta\in C^{*}_{\max}(X)$ be the Laplacian and let $\lambda: C^{*}_{\max}(G(X))\to C^*_u(X)$ be the extension of the regular representation. By assumptions, we know that both $q(\Delta)$ and $\lambda(\Delta)$ have spectral gap, that is
    \[
    \sigma_{C^*_{\max}(\partial G(X))}(q(\Delta))\cup\sigma_{\mathcal B(\ell^2X)}(\lambda(\Delta)) \subseteq \{0\}\cup[c,\infty)
    \]
    for some $c > 0$. Let now $\varphi: [0,\infty)\to [0,1]$ be an arbitrary continuous function with $\supp (\varphi)\subseteq (0,c)$. Consider the element $a = \varphi(\Delta) - \varphi(\Delta)^2 \in C^*_{\max}(G(X))$, which is positive, since $\varphi$ is bounded by $1$. Because of the spectral gap we calculate
    \[
    q(a)=\varphi(q(\Delta))-\varphi(q(\Delta))^2=0.
    \]
    Exactness then implies $a\in \mathcal K(\ell^2 X)$ and thus $\lambda(a)=a$. But similarly we obtain
    \[
    \lambda(a)=\varphi(\lambda(\Delta))-\varphi(\lambda(\Delta))^2=0,
    \]
    so that in fact $a=\lambda(a)=0$. This means $\varphi(\Delta)$ is a projection, so that by the spectral mapping theorem $\varphi(\sigma_{C^*_{\max}(G(X))}(\Delta))=\sigma_{C^*_{\max}(G(X))}(\varphi(\Delta))\subseteq \{0,1\}$. Since $\varphi$ was arbitrary we conclude that $\Delta$ has no spectrum in $(0,c)$ and this finishes the proof.
\end{proof}

Although this proposition justifies the introduction of boundary property (T) in the case of infinite-valued metrics, its usefulness so far is limited by the fact that we can only attest boundary property (T) alongside geometric property (T); we have no example of a space with boundary property (T) that is not also an expander. However, it will be useful for extending the implication $(ii)\Rightarrow (i)$ above to spaces with almost boundary (T). For this we further need the following statement from \cite{Winkel}:

\begin{Lemma}[{see \cite[Proposition 3.1]{Winkel}}]
    Let $X=\bigsqcup_i X_i$ be a box space of graphs with degree bounded uniformly by $d$ such that each $X_i$ is connected and there is $c>0$ with $\sigma_{\mathcal B(\ell^2 X)}(\Delta)\subseteq {0}\cup[c,\infty)$, that is X is an expander. Let $1>\varepsilon>0$, $F_i\subset X_i$, $F=\bigcup_i F_i$ such that $\mu (\overline{F})<\varepsilon$ and $\pi:\mathbb{C}[X]\to \mathcal B(\mathcal H)$ a representation, $v\in \mathcal{H}$ a unit vector, $\eta >0$ such that $\pi(\Delta) v=\eta v$. If $p_F$ is the projection onto the subset $F$ then we have
    \[
    \norm{\pi(p_F) v}^2\leq \frac{2^{\frac{3}{2}}d}{c}\sqrt{\eta}+3\sqrt{\varepsilon}.
    \]
\end{Lemma}

This is slightly more general than \cite[Proposition 3.1 (i)]{Winkel}, but the proof is still essentially the same, so we don't reproduce it here. Note that in \cite{Winkel} the graphs $X_i$ are also required to be connected. 


\begin{Prop} \label{t+exp}
    Let $X=\bigsqcup_i X_i$ be a box space of graphs that is an expander sequence and has almost boundary (T). Assume additionally that each $X_i$ is connected. Then $X$ has geometric property (T).
\end{Prop}

\begin{proof}
    Using Proposition \ref{exp+bdr} it suffices to show spectral gap for the Laplace operator $\Delta\in C^*_{\max}(\partial G(X))$. Consider the short exact sequence
    \begin{equation*}
        0\to C^*_{\max}(\partial G(X)|_{\partial \beta X\setminus Z}) \to C^*_{\max}(\partial G(X))\overset{q}{\to} C^*_{\max}(\partial G(X)|_{Z})\to 0
    \end{equation*}
    of maximal C$^*$-algebras, where $Z\subseteq \pb X$ is a core such that 
    \[
    q(\Delta)\in C^*_{\max}(\partial G(X)|_{Z})
    \]
    has spectral gap, say $\sigma_{C^*_{\max}(\partial G(X)|_{Z})}(q(\Delta))\subseteq \{0\}\cup[c,\infty)$ for $c>0$.

    Assume now that we can find $c>\eta>0$ arbitrarily small lying in the spectrum of $\Delta$ and choose a representation $\pi:C^*_{\max}(\partial G(X))\to \mathcal B(\mathcal H)$ such that $\pi(\Delta)v=\eta v$ for some unit vector $v\in \mathcal H$. Let now $\varphi:[0,\infty)\to[0,1]$ be an arbitrary continuous function with $\supp (\varphi) \subset(0,c)$ and $\varphi(\eta)=1$. Consider the element $\varphi(\Delta)\in C^*_{\max}(\partial G(X))$, so that by spectral gap we have $q(\varphi(\Delta))=\varphi(q(\Delta))=0$ and thus $\varphi(\Delta)\in C^*_{\max}(\partial G(X)|_{\partial \beta X\setminus Z})$ by exactness.

    Take $1>\varepsilon>0$ and choose $f \in C_c(\partial G(X)|_{\partial \beta X\setminus Z})$ with $\norm{\varphi(\Delta)-f}<\varepsilon$. Since $\supp (f)$ is compact and the source and range map are continuous, we know that the set
    \[
    S= r(\supp (\varphi))\cup s(\supp (\varphi))\subset \pb X\setminus Z
    \]
    must also be compact. Recall that the core $Z$ has full measure, so that $\mu_{\mk u}(S)=0$ for any $\mk u\in\partial \beta \mb N$ that we fix. Regularity of $\mu_{\mk u}$ then implies the existence of an open set $U\subseteq \beta X$ with $S\subseteq U$ and $\mu_{\mk u}(U)<\varepsilon$. We can write $U=\bigcup_j U_j$ as union of clopen sets $U_j$, since these form a basis of the topology of $\beta X$. As $S$ is compact the union of finitely many of the $U_j$, which we will denote by $\overline{F}$ for some $F\subseteq X$, will cover $S$ and still be clopen. In particular we have $S\subseteq \overline{F}\subseteq U$ and thus $\mu_{\mk u}(\overline{F})<\varepsilon$. Consider $p_F$, the projection onto $F$, as an element of $C^*_{\max}(\partial G(X))$, then since $S\subseteq \overline{F}$ we have $f p_F=p_F f=f$. The preceding lemma together with $\norm{f}\leq 2$ then allows us to obtain the following estimate:
    \begin{align*}
        \norm{\pi(f)v}^2 &= \norm{\pi(f)\pi(p_F)v}^2\leq 4\norm{\pi(p_F)v}^2 \leq C_1\sqrt{\eta}+C_2\sqrt{\varepsilon}
    \end{align*}
    for appropriate constants $C_1,C_2>0$. Together with
    \[
    \norm{\pi(\varphi(\Delta))v}\leq \norm{\varphi(\Delta)-f}+\norm{\pi(f)v}\leq \varepsilon+\norm{\pi(f)v}
    \]
    letting $\varepsilon\to 0$ gives us $\norm{\pi(\varphi(\Delta))v}^2\leq C_1\sqrt{\eta}$. But because of $\varphi(\eta)=1$ we also have 
    \[
    \pi(\varphi(\Delta))v=\varphi(\pi(\Delta))v=\varphi(\eta)v=v
    \]
    and thus $1=\norm{v}^2\leq C_1\sqrt{\eta}$. This is absurd since $\eta$ was to be chosen arbitrarily small. So $\Delta$ must have had a spectral gap in the first place.
\end{proof}

\begin{Rem}
    There are two things worth observing about the proof of Proposition \ref{t+exp}. First, it did not make use of the definition of the core, only that it is a closed and invariant subset of $\pb X$ of full measure. However, the precise structure of the core will play a role in the next to sections. Second, while the spectral gap need not be preserved like in the previous proposition, there is still a lower bound to it. Indeed, we arrive at a contradiction for $\eta < \frac{1}{C_1^2}$, so that with the constants $d$ and $c$ from before the spectral gap is at least $\frac{c^2}{144d^2}$.
\end{Rem}

As a corollary we obtain an alternative proof of a result from \cite{Winkel}:

\begin{Cor}[{\cite[Theorem B]{Winkel}}]
    Let $X$ and $X'$ be approximately isomorphic box spaces of connected graphs such that $X$ has geometric property (T) and $X'$ is an expander. Then $X'$ has geometric property (T).
\end{Cor}

\begin{Ex}\label{examplenonconnectedexp+almostboundaryT}
    The requirement that $X_i$ be connected is indeed necessary: take a box space $X$ with geometric property (T) and an expander box space $Y$ without geometric property (T) such that $|Y_i|/|X_i|\to 0$. Form graphs $X_i'$ that are the disjoint unions of the graphs $X_i$ and $Y_i$, then the resulting box space $X'$ is approximately isomorphic to $X$ and $G(X')$ has $\beta X$ and $\beta Y$ as invariant sets. Then $X'$ cannot have geometric property (T) since spectral gap of $\Delta\in C^*_{\max}(G(X'))$ would pass to $C^*_{\max}(G(Y))\cong C^*_{\max}(G(X')|_{\beta Y})$. 
\end{Ex}

As is turns out this counterexample is already the worst thing that can happen. That is, by neglecting the bad parts we will always obtain a space with geometric property (T), as the next result will tell us.

\begin{Thm}\label{exp+bdrDiffuse}
    Let $X=\bigsqcup_i X_i$ be a box space of graphs that is an expander and has almost boundary (T). Then there is a box space $X'$ that is approximately isomorphic to $X$ and has geometric property (T).
\end{Thm}

\begin{proof}
    The idea of proof is to sort out all connected components of the $X_i$ that are not properly witnessed by the sets $Y_i$ which generate the core for almost boundary (T) of $X$. The resulting space can then be arranged to make use of the previous proposition in such a way that boundary property (T) carries over.

    We denote by $X_{i,k}$ for $k=1,\ldots,K_i$ the connected components of $X_i$. Take $\varepsilon_i\to 0$ with $|Y_i|/|X_i|=1-\varepsilon_i$ and write $W_{i,k}=X_{i,k}\setminus Y_i$. By rearranging we can assume $|W_{i,k}|\geq\sqrt{\varepsilon_i}|X_{i,k}|$ for $k=1,\ldots,L_i$ and $|W_{i,k}|<\sqrt{\varepsilon_i}|X_{i,k}|$ for $k=L_i+1,\ldots,K_i$. We estimate
    \[
    \sum_{K=1}^{L_i} |X_{i,k}|\leq\frac{1}{\sqrt{\varepsilon_i}}\sum_{k=1}^{L_i}|W_{i,k}|\leq \frac{1}{\sqrt{\varepsilon_i}}\sum_{k=1}^{K_i}|W_{i,k}|=\frac{1}{\sqrt{\varepsilon_i}}|X_i\setminus Y_i|=\sqrt{\varepsilon_i}|X_i|
    \]
    and conclude that for $X_i'=\bigcup_{k=L_i+1}^{K_i}X_{i,k}$ we obtain a box space $X'=\bigsqcup_i X_i'$ that is approximately isomorphic to $X$.

    Define a box space $\tilde{X}=\bigsqcup_j \tilde{X}_j$ where the $\tilde{X}_j$ are the connected components $X_{i,k}$ of $X'$ in lexicographic order (with $i$ taking precedence over $k$ of course). Let $\tilde{Y}_j=Y\cap \tilde{X}_j$ for $Y=\bigcup_i Y_i$ so that by the above construction we have $|\tilde{Y}_j|/|\tilde{X}_j|>1-\sqrt{\varepsilon_{j(i)}}\to 1$. By definition $\beta X'=\beta \tilde{X}$ and thus the sets $\partial E_R$ must also coincide, yielding $\partial G(X')=\partial G(\tilde{X})$. So Proposition \ref{t+exp} applies to $\tilde{X}$. In particular, $X'$ must have boundary property (T). As $X'$ is obviously an expander it must also have geometric property (T) by Proposition \ref{exp+bdr}.
\end{proof}

\section{Measured property (T)}\label{section: measured T}

We introduce a definition of measured property (T) for p.m.p. groupoids inspired by \cite{dA-W} which is equivalent to the one given in \cite{A-D} for standard Borel p.m.p. groupoids but also works in the non-standard case and is therefore well-suited for our purposes. We will recall some definitions from \cite{A-D} and \cite{dA-W} and refer to these articles for further background on property (T) for groupoids.

\begin{Def}
Let $G$ be a discrete p.m.p. measured groupoid. A $\ast$-representation $\pi \colon \mb C G\to \mb B(\mc H)$ is called normal if its restriction to $L^\infty(G^{(0)})$ is normal (that is, weak$^\ast$-continuous). A vector $\xi\in \mc H$ is called invariant or constant if for every $f\in\mb CG$ we have $\pi(f)\xi = \pi(\Psi(f))\xi$. The subspace of constant vectors in $\mc H$ is denoted by $\mc H^\pi$, and we denote its orthogonal complement by $\mc H_\pi$.
\end{Def}

\begin{Def}
A discrete p.m.p. measured groupoid $(G,\mu)$ has (measured) property (T) if there exists $c>0$ and finitely many bisections $\phi_1,\dots,\phi_n$ which generate $G$ such that for every normal $\ast$-representation $\pi\colon \mb C G\to \mb B(\mc H)$ and for all $\xi\in \mc H_\pi$ (orthogonal complement of invariant vectors) we have
\[
\norm{(\phi_i - \Psi(\phi_i))\xi}\geq c \norm{\xi}
\]
for at least one $i\in\{1,\dots,n\}$.
\end{Def}

\begin{Prop}
Let $G$ be a standard Borel $r$-discrete p.m.p. measured groupoid. Then $G$ has property (T) in the sense of \cite{A-D} if and only if it has property (T) in the above sense.
\end{Prop}
\begin{proof}
($\implies$): we prove the contrapositive following the proof in \cite[Theorem 12.1.7]{brown-ozawa} for groups. Let $E_1\subset E_2\subset \cdots\subset G$ be an exhaustion of $G$ by generating subsets of finite measure; we identify each $E_k$ with a finite set of disjoint bisections. If measured property (T) fails, then 
every $k$ there exists a normal $\ast$-representation $\pi_k\colon \mb C G\to \mb B(\mc H_k)$ and $\xi_k\in \mc H_k$ such that
\[
\delta_k\coloneqq \sup_{\phi\in E_k}\norm{(\phi - \Psi(\phi))\xi_k}< 4^{-k} \norm{\xi_k-P_k\xi_k},
\]
where $P_k$ is the orthogonal projection onto the subspace of invariant vectors. We set $\mc H\coloneqq \bigoplus_k \mc H_k$ and $\pi\coloneqq \bigoplus_k \pi_k$. Then $\pi$ is a normal representation of $\mb C G$ on $\mc H$. We now define a cocycle
\[
\sigma(\gamma)\coloneqq \left(\frac{ \xi_k\circ r(\gamma) - L_{\pi_k}(\gamma)\xi_k\circ s(\gamma)}{2^k\delta_k}\right)_{k\in \mb N}
\]
with values in $\mc H$ which is not inner on $G$ since $\norm{\sigma(\phi)}\geq 2^k$ for $\phi\in E_k$. 

($\impliedby$): Suppose that $\eta_k$ is a sequence of almost invariant unit vectors in $\mc H$ and let $\xi_k\coloneqq \eta_k - P\eta_k\in \mc H_\pi$, where $P$ is the orthogonal projection onto the subspace of invariant vectors. Then for every $i = 1,\dots,n$ we have $\norm{(\phi_i - \Psi(\phi_i))\xi_k}\to 0$ which by the definition of measured property (T) implies $\norm{\xi_k}\to 0$. But then $P\eta_k\neq 0$ for $k$ big enough, providing the desired invariant vector.
\end{proof}


\begin{Lemma}\label{lem:kernel-const}
Let $(G,\mu)$ be a discrete p.m.p. measured groupoid and $\phi_1,\dots,\phi_n$ finitely many bisections which generate $G$. Then for every normal representation $\pi$ of $\mb CG$ the kernel of the Laplace operator $\pi(\Delta)$ in that representation is exactly the space of constant vectors.
\end{Lemma}
\begin{proof}
The proof is a straightforward adaptation of \cite[Lemma 7.6]{dA-W}, even with some simplifications in the measurable case, so we only sketch the argument.

If $\xi$ is a constant vector, then clearly $\pi(\Delta)\xi = 0$. Conversely, if $\pi(\Delta)\xi = 0$, then $\pi(\phi_i)\xi = \Psi(\phi_i)\xi$ for all $i=1,\dots,n$. As these bisections are generating, this implies that $\pi(\phi)\xi = \Psi(\phi)\xi$ for every bisection $\phi$ of $G$, which in turn implies that $\xi$ is constant.
\end{proof} 

\begin{Prop}\label{proposition: measured T characterization}\quad 
\begin{enumerate}
\item A compactly generated étale topological groupoid $G$ has (topological) property (T) if and only if $\Delta$ has spectral gap in every representation of $G$.
\item A discrete p.m.p. measured groupoid $(G,\mu)$ has (measured) property (T) if and only if the Laplacian $\Delta$ has spectral gap in every normal representation of $\mb CG$.
\end{enumerate}
\end{Prop}
\begin{proof}
The first part follows from \cite[7.5--7.7]{dA-W} and the observation that \cite[Lemma 7.5]{dA-W} actually characterizes topological property (T).

The second part can be deduced in \cite[Lemma 7.7]{dA-W}, using Lemma \ref{lem:kernel-const} above. Indeed, for every normal representation $\pi\colon \mb CG\to B(\mc H)$ and every $\xi\in \mc H$ we have
\[
\ip{\pi(\Delta) \xi,\xi} = \sum_{i=1}^n \norm{(\phi_i - \Psi(\phi_i))\xi}^2.
\]
Since the kernel of $\pi(\Delta)$ is exactly the space of constant vectors, the existence of a spectral gap $c>0$ for  $\pi(\Delta)$ is equivalent to the existence such that for every $\xi\in \mc H_\pi$ we have $\ip{\pi(\Delta) \xi,\xi}\geq c^2\norm{\xi}^2$, which is in turn equivalent to the condition in the definition of measured property (T).
\end{proof}

The following proposition extends \cite[Theorem 5.12]{dA-W} to the non-ergodic and non-standard case.

\begin{Prop}
Let $G$ be a compactly generated étale topological groupoid with (topological) property (T) and $\mu$ an invariant measure on $G^{(0)}$. Then the discrete p.m.p. measured groupoid $(G,\mu)$ has (measured) property (T). 
\end{Prop}
\begin{proof}
By \cite[7.5--7.7]{dA-W} the Laplace operator for $G$ has spectral gap in every representation of $G$. Since every normal representation of $\mb CG$ is a representation of $G$, the result follows.
\end{proof}

\begin{Cor}\label{corollary: almost t implies measured T}
If a box space of graphs $X$ has almost boundary (T), then for every ultrafilter $\mk u\in \partial\beta\mb N$ the discrete p.m.p. measured groupoid $(\partial G(X),\mu_{\mk u})$ has measured property (T).
\end{Cor}

\section{Spectral gap in the ultraproduct and expanders}

It appears now very desirable to show that almost boundary (T) or measured (T) for a box space guarantees being approximately isomorphic to an expander. This is what we will show in this section by adapting the method of Kun \cite{kun}. We will keep our text self-contained since we will need appropriate modifications of the statements and proofs, referring to analogous results from his paper.

Let us first recapitulate the main idea of Kun's work: using  spectral gap coming from property (T) one can find a constant $C>0$ and partition every graph in the sequence into components such that their subsets have expansion at least $C$ and such that there are only a few edges between these components. To produce these components, one starts with a characteristic function $\chi_T$ of a subset with expansion less than $C$ and repeatedly applies the Markov operator; spectral gap then guarantees that some level set $U_t = \{M^K\chi_T(x) > t\}$ will be a desired component. At the end, one runs a rewiring argument to obtain a disjoint union of graphs that make up our expander. 

Our main result will only use spectral gap in the ultraproduct which is a weaker condition, so we will be able to apply the above above construction only for subsets $T$ of sufficiently big measure. This requires a modification for the case of ``small'' subsets $T$ in Propositions \ref{UforT} and \ref{almostexpanderdecomposition} below.

\begin{Def}
    Let $X=\bigsqcup_iX_i$ be a box space and $\mk u\in \partial\beta \mb N$ some non-principal ultrafilter. We say that $X$ has \emph{spectral gap $c>0$ in the ultraproduct} along $\mk u$, if the Laplace operator $\Delta$ has spectral gap $c>0$ as an operator on $\leftsub{{i\to\mk u}}{\prod} L^2(X_i,\mu_i)$. Furthermore we say that $X$ is an expander along $\mk u$, if there is $c>0$ and $I\in \mk u$ such that $\Delta\in \mathcal B(L^2(X_i,\mu_i))$ has spectral gap $c$ for every $i\in I$.
\end{Def}

\begin{Rem}
    From Proposition \ref{proposition: measured T characterization} it is clear that if the measured groupoid $(\partial G(X),\mu_{\mk u})$ has measured property (T), then $X$ must have spectral gap in the ultraproduct along $\mk u$. In particular, by Corollary \ref{corollary: almost t implies measured T}, this holds true for every $\mk u\in \partial \beta \mathbb N$ if $X$ has almost boundary (T).
\end{Rem}

Obviously, an if $X$ is an expander along $\mk u\in\partial\beta \mb N$, then it has spectral gap in the ultraproduct along $\mk u$. Our main result in this section is the converse:

\begin{Thm}\label{theorem: spectral gap implies almost exp}
    Let $X$ be a box space that has spectral gap in the ultraproduct along some $\mk u\in \partial\beta \mb N$. Then there is a box space $X'$, that is approximately isomorphic to $X$ and that is an expander along $\mk u$.
\end{Thm}

By standard arguments for ultrafilters this implies the following:

\begin{Cor}\label{corollary:almost T is appr expander}
    A box space $X$ has spectral gap in the ultraproduct along every ultrafilter $\mk u\in \partial\beta \mb N$ if and only if $X$ is approximately isomorphic to an expander. In particular, if $X$ has almost boundary (T), then it is approximately isomorphic to an expander.
\end{Cor}

Using the results of the previous sections we arrive at the theorem we set out to obtain.

\begin{Thm}\label{theorem: almost T implies appr eq to geometric T}
    Let $X$ be a box space with almost boundary (T). Then there is a box space $X'$ that is approximately isomorphic to $X$ and that has geometric property (T).
\end{Thm}

\begin{proof}
    From Corollary \ref{corollary:almost T is appr expander} $X$ is approximately isomorphic to an expander $X''$, which has almost boundary (T) by Proposition \ref{almTinv} . Theorem \ref{exp+bdrDiffuse} then implies the existence of geometric property (T) box space $X'$, that is approximately isomorphic to $X''$ and thus to $X$.
\end{proof}

So let us turn to the proof of Theorem \ref{theorem: spectral gap implies almost exp}. We follow the construction of Kun in \cite{kun} and adjust the argumentation to hold in this more general setting. Since the graphs in our sequences are not $d$-regular but only of degree bounded by $d$, we need to define the Markov operator in an appropriate way:

\begin{Def}
    For a box space $X$ we define $M=\mathbf 1 -\frac{1}{2d}\Delta\in \mathbb C[X]$.
\end{Def}

Note that we scaled $\Delta$ with the factor $\frac{1}{2d}$ to guarantee that $M$ has spectral gap at $1$ if and only if $\Delta$ has spectral gap at $0$ in any representation of $\mathbb C[X]$. This works out since we can estimate the maximal norm of $\Delta$ by
\[
\norm{\Delta}\leq \sup_{x\in \beta X}\sum_{y\in G(X)_x}|\Delta_{xy}|\leq \sup_{x\in \beta X}(|\Delta_{xx}|+d\cdot 1)\leq 2d,
\]
see e.g. \cite[p. 205]{brown-ozawa}. We collect some basic properties of $M$ in a lemma.

\begin{Lemma}[{compare with \cite[Lemma 7]{kun}}]
    Take $M=\mathbf 1 -\frac{1}{2d}\Delta\in \mathbb C[X]$ as above and let $\pi:\mathbb C[X]\to\mathcal B(\mathcal H)$ be a representation of $\mathbb C[X]$. Then $\norm{M}\leq 1$ and $\ker \pi(\Delta)=\mathrm{fix}\,\pi(M)$. If moreover $\pi(\Delta)$ has spectral gap $c>0$ at $0$, then $\pi(M)$ has spectral gap $c_M=\frac{c}{2d}$ at $1$ and we have
    \[
    \norm{\pi(M)^{k+1}f-\pi(M)^kf}\leq \left(1-c_M\right)^k\norm{\pi(M)f-f}
    \]
    for every $f\in \mathcal H$ and $k\in \mathbb N$.
\end{Lemma}

\begin{proof}
    The definition of $M$ immediately yields $\ker \pi(\Delta)=\mathrm{fix}\,\pi(M)$. From the spectral mapping theorem we have that $\sigma_{\mathcal B(\ell^2X)}(M)=1-\frac{1}{2d}\sigma_{\mathcal B(\ell^2X)}(\Delta)\subseteq[0,1]$, since $\sigma_{\mathcal B(\ell^2X)}(\Delta)\subseteq[0,2d]$. Since $M$ is self-adjoint, we have $\norm{M}=\sup_{\lambda\in\sigma_{\mathcal B(\ell^2X)}(M)}|\lambda|\leq 1$. The statement about the spectral gap follows similarly.
    To show the inequality denote $c_M=\frac{c}{2d}$, then for every $f \perp \fix\pi(M)$ we have 
    \[
    \norm{\pi(M)f}\leq(1-c_M)\norm{f}.
    \]
    Since $\scp{\pi(M)f-f,h}=\scp{f,\pi(M)h}-\scp{f,h}=0$ for every $h\in \fix\pi(M)$, we have $\pi(M)f-f\perp \fix\pi(M)$ for every $f\in\mathcal H$, which yields
    \[
    \norm{\pi(M)^2f-\pi(M)f}\leq(1-c_M)\norm{\pi(M)f-f}
    \]
    and thus the assertion follows by induction.
\end{proof}

The next lemma will allow us to relate spectral gap in the limit object to certain subsets of the $X_i$.


\begin{Lemma}[{compare with \cite[Lemma 10]{kun}}]\label{ultralemma}\quad
 Let $X=\bigsqcup_i X_i$ have 
 spectral gap $c>0$ in the ultraproduct along $\mk u$. Then for each $\varepsilon>0$, each $\delta > 0$ and $k\in \mathbb N$ there is $I\in \mk u$ such that for all $i\in  I$ and for every $T\subseteq X_i$ with $\mu_i(T)\geq \delta$ we have 
    \beq\label{eq:ultralemma}
    \norm{M^{k+1}\chi_T-M^k\chi_T}\leq ((1-c_M)^k+\varepsilon)\norm{M\chi_T-\chi_T},
    \eeq
    where $c_M=c/2d$.

Furthermore, if $X$ has almost boundary (T) with respect to $Z=\mathrm{core}(Y)$, then for each $\varepsilon>0$ and $k\in \mathbb N$ here exists $R > 0$ and $N\in \mb N$ such that for all $i\geq N$ and for every $T\subseteq X_i$ with $T\subseteq X_i\setminus B_R(X_i\setminus Y_i)$ the inequality \eqref{eq:ultralemma} holds.
\end{Lemma}

\begin{proof}
    We prove the second statement first. To achieve a contradiction let there be $\varepsilon>0$ and $k\in \mathbb N$ such that for any $R\in \mathbb N$ we can find $i_R\in \mathbb N$ and $T_R\subseteq X_{i_{R}}\setminus B_R(X_{i_{R}}\setminus Y_{i_{R}})$ with
    \[
    \norm{M^{k+1}\chi_{T_{R}}-M^k\chi_{T_{R}}} > ((1-c')^k+\varepsilon)\norm{M\chi_{T_{R}}-\chi_{T_{R}}}.
    \]
    Clearly this implies $i_R\to\infty$ as $R\to\infty$. We define a measure $\nu_R$ on $X_{i_{R}}$ by $\nu_R(A)=\frac{1}{|T_R|}|A|$. Take some non-principal ultrafilter $\omega \in \pb \mathbb N$ and let $\nu =\lim_{R\to\omega}\nu_R$ and declare $T= \pb\left(\bigcup_R T_R \right)\subseteq\pb X$. By construction we then have $T\subseteq Z$ and $\nu(T)=1$ so that $\chi_T$ is a unit vector in $L^2(Z,\nu)$. Denote by $\pi:C_c(G(X)|_Z)\to \mathcal B(L^2(Z,\nu))$ the standard representation (which is just the fiber-wise application of an element in $C_c(G(X)|_Z)$). By the preceding lemma we thus have
    \[
    \norm{\pi(M)^{k+1}f-\pi(M)^kf}\leq \left(1-c'\right)^k\norm{\pi(M)f-f}.
    \]
    But clearly $\pi(M)$ is also the respective limit operator on the ultraproduct $\prod_{R\to \omega}L^2(X_{i_{R}},\nu_R)$ and hence satisfies 
    \[
    \norm{\pi(M)^{k+1}\chi_T-\pi(M)^k\chi_T} \geq ((1-c')^k+\varepsilon)\norm{\pi(M)\chi_T-\chi_T},
    \]
    which is the desired contradiction.

    The first statement follows similarly by observing that in this case the vector $\chi_T$ belongs to the ultraproduct representation (in fact, $\chi_T$ has norm at least $\sqrt \delta$ there) so that we can use the spectral gap in this representation.
\end{proof}

We will also need an estimate for the size of the boundary of level sets.

\begin{Lemma}[{compare with \cite[Lemma 8]{kun}}]\label{boundary of U lemma}
    Let $G$ be a finite graph with degree bounded by $d$, let $0<a<b<1$ and $f\in\ell^2G$. Then there is $t\in(a,b)$ such that for $U=\{x\in G\,:\,f(x)>t\}$ we have 
    \[
    |\partial U|^2\leq\frac{4d^2}{a^2(b-a)^2}\norm{Mf-f}\norm{f}^3.
    \]
\end{Lemma}

\begin{proof}
    Recall the co-area formula
    \[
    \int^\infty_0|\partial\{f>t\}|\dt=\frac{1}{2}\sum_{x\in G,\,x\sim y}|f(x)-f(y)|.
    \]
    From this, the Cauchy-Schwarz inequality and Markov's inequality we estimate for every $t\in (a,b)$
    \begin{align*}
        \left(\frac{1}{b-a}\int^b_a|\partial \{f>t\}|\dt\right)^2
        &\leq \left(\sum_{f(x)\geq a}\sum_{y\sim x}\frac{|f(x)-f(y)|}{b-a}\right)^2\\
        &\leq \sum_{f(x)\geq a}\sum_{y\sim x}1^2\sum_{x\in G}\sum_{y\sim x}\left(\frac{|f(x)-f(y)|}{b-a}\right)^2\\
        &\leq d|\{f\geq a\}|\frac{1}{(b-a)^2}2\langle \Delta f,f\rangle\\
        &\leq d\frac{\|f\|^2}{a^2}\frac{1}{(b-a)^2}2\|\Delta f\|\|f\|\\
        &=\frac{4d^2}{a^2(b-a)^2}\|Mf-f\|\|f\|^3.
    \end{align*}
    Since this is true on average there must be some $t\in (a,b)$ such that $U=\{f>t\}$ has the desired properties.
\end{proof}

Using these lemmas we can formulate the result that was suggested at the beginning of this section: if we have a ``good" subset $T$ with small boundary we can replace it by a similar set $U$ with even smaller boundary.

\begin{Prop}[{compare with \cite[Corollary 9 and Proposition 11]{kun}}]\label{UforT}
    Let $X=\bigsqcup_i X_i$ have spectral gap $c > 0$ in the ultraproduct along $\mk u$. Then there is a constant $C>0$ such that for all $\alpha >0$ there is $K > 0$ such that for every $\delta > 0$ there is $I\in \mk u$ so that for all $i\in I$ and every $T\subseteq X_i$ with $\mu_i(T)\geq \delta$ and $|\partial T|<C|T|$ there exists $U\subseteq B_K(T)$ with $|U\triangle T|<|T|/4$ and $|\partial U|<\alpha |U|$.
\end{Prop}

\begin{proof}
    Take $C=\frac{c_M^2}{72}$. Let $\alpha>0$ and take $\varepsilon>0$ and $K\in \mathbb N$ such that $(1-c_M)^K\leq \frac{\alpha^2}{5184d^2}$ and $\varepsilon=\min\left\{\frac{\alpha^2}{5184d^2},\frac{1}{Kc_M}\right\}$. Now given $\delta > 0$ there is a $I\in \mk u$ such that the conclusion of Lemma \ref{ultralemma} holds for each $1\leq k\leq K$.

    Now let $i\in I$ and $T\subseteq X_i$ with $\mu_i(T)\geq \delta$ and $|\partial T|<C|T|$. Write $f=M^K\chi_T$ as well as $a=\frac{1}{3}$ and $b=\frac{2}{3}$. We use Lemma \ref{boundary of U lemma} to obtain a set $U=\{f>t\}$ with $t\in(a,b)$. Since $f(x)>\frac{1}{3}$ for $x\in U$ and $f(x)<\frac{2}{3}$ for $x\notin U$ we have $|f(x)-\chi_T(x)|>\frac{1}{3}$ for each $x\in U\triangle T$. This implies   
    \[
    |U\triangle T|\leq\sum_{x\in U\triangle T}9|f(x)-\chi_T(x)|^2\leq 9\norm{M^K\chi_T-\chi_T}^2.
    \]
    Furthermore we estimate
    \begin{align*}
    \norm{M^K\chi_T-\chi_T}&\leq\sum_{k=0}^{K-1}\norm{M^{k+1}\chi_T-M^k\chi_T}\leq\sum_{k=0}^{K-1}((1-c_M)^k+\varepsilon)\norm{M\chi_T-\chi_T}\\
    &\leq \left(\sum_{k=0}^\infty(1-c_M)^k+K\varepsilon\right)\norm{M\chi_T-\chi_T}\\
    &=\left(\frac{1}{c_M}+K\varepsilon\right)\norm{M\chi_T-\chi_T}
    \end{align*}
    and this gives us
    \begin{align*}
    |U\triangle T| &\leq 9\left(\frac{1}{c_M}+K\varepsilon\right)^2\norm{M\chi_T-\chi_T}^2\leq 9\left(\frac{2}{c_M}\right)^2\frac{1}{(2d)^2}\norm{\Delta\chi_T}^2\\
    &= \frac{9}{(c_Md)^2}\sum_{x\in \supp(\partial T)}|\Delta\chi_T(x)|^2\leq \frac{9}{(c_Md)^2} |\supp(\partial T)|d^2\\
    &\leq \frac{18}{c_M^2}|\partial T| < \frac{18}{c_M^2}C|T|\leq \frac{|T|}{4}
    \end{align*}
    by the choice of $\varepsilon$ and $C$.

    Finally, since
    \[
    |T|=|T\cap U|+|T\setminus U|\leq |U|+|U\triangle T|\leq |U|+\frac{|T|}{2}
    \]
    we have $|T|\leq 2|U|$ and so for the boundary of $U$ we obtain from Lemma \ref{ultralemma}
    \begin{align*}
    |\partial U|^2 &\leq \frac{4d^2}{\left(\frac{1}{3}\right)^2\left(\frac{2}{3}-\frac{1}{3}\right)^2}\norm{M^{K+1}\chi_T-M^K\chi_T}\norm{M^K\chi_T}^3\\
    &\leq 324d^2((1-c_M)^K+\varepsilon)\norm{M\chi_T-\chi_T}\norm{\chi_T}^3\\
    &\leq 324d^2((1-c_M)^K+\varepsilon)2|T|^{\frac{1}{2}}|T|^{\frac{3}{2}}\\
    &= 648d^2((1-c_M)^K+\varepsilon)|T|^2\\
    &\leq 648d^2((1-c_M)^K+\varepsilon)4|U|^2\\
    &\leq \alpha^2|U|^2
    \end{align*}
    by the choice of $K$ and $\varepsilon$.
\end{proof}

We will use this result to decompose our box space into subsets with a lot of ``inner" expansion and only few edges between them. The crucial step is to deal with those subsets $T\subseteq X_i$ for which $\mu_i(T)<\delta$. This is done by considering such a $T$ as ``bad", if it cannot be replaced by an adequate $U$, and showing that those ``bad" subsets in total have measure less than $\delta$.


\begin{Prop}[{compare with \cite[Theorem 3]{kun}}]\label{almostexpanderdecomposition}
    Let $X=\bigsqcup_iX_i$ have spectral gap in the ultraproduct along $\mk u$. Then there is a constant $C>0$ such that for every $\alpha>0$ exists $I\in \mk u$ such that for all $i\in I$ the following holds true: there is a partition $X_i=\bigsqcup_{\ell=0}^{L_i}P_i^\ell$ such that $|P_i^{0}|<\alpha |X_i|$ and for every $\ell = 1,\dots,L_i$ we have $|\partial P_i^\ell|<\alpha |P_i^\ell|$ as well as
    \[
    T\subseteq P_i^\ell\text{ and } |T|\leq |P_i^\ell|/2 \implies |\partial T|\geq C|T|.
    \]
\end{Prop}

\begin{proof} 
    We adopt the notation of Proposition \ref{UforT}. So take $C$ as before and let $\alpha>0$ be given (w.l.o.g. we may assume $\alpha<1$). We apply the previous proposition to $\alpha'=\frac{\alpha^4}{10^4}$ and find a $K > 0$ such that for $\delta\coloneqq \frac{1}{d^{3K+1}}\frac{\alpha^4}{10^4}$ the conclusion of the proposition holds.
    
    Let's call a subset $T\subseteq X_i$ \emph{good} if there exists $U\subset X_i$ with $|U\triangle T|<|T|/2$ and $|\partial U|\leq \frac{\alpha^2}{10^2}|U|$.  We inductively construct a decomposition of $X_i$ for $i\in I$ in the following way. Initialize $P_i^{0}\coloneqq \varnothing$.  Suppose now $P_i^{0},\ldots,P_i^{L}$ to be defined already. If $|\partial T|\geq C|T|$ for every $T\subseteq X_i\setminus(\bigsqcup_{\ell}^LP_i^{\ell})$, we set $P_i^{L+1}=X_i\setminus(\bigsqcup_{\ell}^LP_i^{\ell})$ and $L_i=L+1$. If not, take some $T=T^{L+1}_i\subseteq X_i\setminus(\bigsqcup_{\ell}^LP_i^{\ell})$ such that $|\partial T|<C|T|$ and $|T|$ is minimal.

    Now there are two cases. If $T$ is not good, we add $B_{3K}(T)\setminus(\bigsqcup_{\ell}^LP_i^{\ell})$ to $P_i^{0}$ and continue the recursion. If $T$ is good, we get a subset $U=U^{L+1}_i\subseteq X_i$ such that $|U\triangle T|\leq|T|/2$ and $|\partial U|\leq \frac{\alpha^2}{10^2}|U|$. In particular we have
    \begin{align*}
        |U| &= |U\setminus T|+|U\cap T|\leq |U\triangle T|+|T|\\
        &\leq\frac{1}{2}|T|+|T|=\frac{3}{2}|T|
    \end{align*}
    and thus $|U|<2|T|$ as well as
    \begin{align*}
        |U\cap T| &= |T\cup U|-|U\triangle T|\geq|T|-\frac{1}{2}|T|\\
        &= \frac{1}{2}|T| \geq \frac{1}{3}|U|.
    \end{align*}
    We define $P^{L+1}_i=U_i^{L+1}\setminus(\bigsqcup_{\ell}^LP_i^{\ell})$. If now $S\subseteq P^{L+1}_i$ with $2|S|\leq|P^{L+1}_i|$, then we have
    \[
    |S|\leq\frac{1}{2}|P^{L+1}_i|\leq\frac{1}{2}|U|<|T|
    \]
    and thus follows $|\partial S|\geq C|S|$ since we chose $T$ to be minimal. Note that because of $|P^{L+1}_i|\geq|U\cap T|\geq\frac{1}{2}|T|>0$ we know that $P^{L+1}_i$ is not empty, so that the step is well defined.

    We now claim that $|P_i^{0}|\leq \frac{1}{d}\frac{\alpha^2}{10^2}|X_i|$. Indeed, let $W_1,\dots,W_m$ be the bad subsets obtained in the above process and let $W = \bigsqcup W_j$. We claim that $\mu_i(W) < \delta$. Supposing the contrary, we use the previous Proposition to find $V\subseteq B_K(W)$ with $|V\triangle W|<|W|/4$ and $|\partial V|<\frac{\alpha^4}{10^4}|V|$. Since the bad subsets $W_1,\dots,W_m$ were $3K$-separated by construction, $V$ decomposes as $V = \bigsqcup_j V_j$ such that $V_j\subseteq B_K(W_j)$ are $K$-separated. We thus have
    \[
        \sum_j |\partial V_j| = |\partial V| <  \frac{\alpha^4}{10^4}|V| = \frac{\alpha^4}{10^4} \sum_j |V_j|
    \]
    and 
    \[
        \sum_j |W_j \triangle V_j| = |W\triangle V| < \frac{|W|}{4} = \sum_j \frac{|W_j|}{4}.
    \]
    The Markov inequality implies that
    \[
    \sum_{j: |W_j \triangle V_j| \leq |W_j|/2}|W_j| \geq |W|/2\geq |V|/4,
    \]
    and consequently
    \[
    \sum_{j: |W_j \triangle V_j| \leq |W_j|/2}|V_j| \geq |V|/8,
    \]
    and 
    \[
    \sum_{j: |\partial V_j| < \frac{\alpha^2}{10^2}|V_j|} |V_j| \geq \left(1-\frac{\alpha^2}{10^2}\right)|V|.
    \]
    Therefore there exists $j$ such that $|W_j \triangle V_j| \leq |W_j|/2$ and $|\partial V_j| < \frac{\alpha^2}{10^2} |V_j|$. But this means that $W_j$ is good, a contradiction. We thus have $\mu_i(W)<\delta$ and consequently $|P_i^{0}|\leq |B_{3K}(W)|\leq d^{3K+1}|W|<\frac{1}{d}\frac{\alpha^2}{10^2}|X_i|$.

    It remains to show the assertion about the boundaries. By construction we have that each edge between some $P_i^\ell$ and $P_i^k$ for $\ell> k>0$ is in particular an edge between $P_i^\ell$ and $U_i^k$ so that we have $\sum_{\ell=1}^{L_i}|\partial P_i^\ell|\leq2\sum_{\ell=1}^{L_i-1}|\partial U_i^\ell|$. From this we have
    \begin{align*}
        \sum_{\ell=0}^{L_i}|\partial P_i^\ell| &\leq d|P_i^{0}|+2\sum_{\ell=1}^{L_i-1}|\partial U_i^\ell|\leq  \frac{\alpha^2}{10^2}|X_i|+2\sum_{\ell=1}^{L_i-1}\frac{\alpha^2}{10^2}| U_i^\ell|\\
        &\leq \frac{\alpha^2}{10^2}|X_i|+6\sum_{\ell=1}^{L_i-1}\frac{\alpha^2}{10^2}| P_i^\ell|\leq 7\frac{\alpha^2}{10^2}|X_i|\leq\frac{\alpha^2}{2}|X_i|.
    \end{align*}
    Without loss of generality we assume now that $|\partial P_i^\ell|\geq \alpha|P_i^\ell|$ for $\ell=1,\dots,L_i'$ and $|\partial P_i^\ell|< \alpha|P_i^\ell|$ for $\ell>L_i'$. Then
    \begin{align*}
        \sum_{\ell=1}^{L'_i}| P_i^\ell| \leq \frac{1}{\alpha}\sum_{\ell=1}^{L'_i}|\partial P_i^l| 
        \leq \frac{1}{\alpha} \sum_{l=0}^{L_i}|\partial P_i^\ell|
        \leq \frac{\alpha^2}{2\alpha}|X_i|=\frac{\alpha}{2} |X_i|.
    \end{align*}
    Thus we have
    \[
    \sum_{\ell=0}^{L'_i}| P_i^\ell| = |P_i^{0}|+\sum_{\ell=1}^{L'_i}| P_i^\ell| \leq \frac{\alpha}{2}|X_i|+\frac{\alpha}{2}|X_i|=\alpha|X_i|.
    \]
    If we now declare $P_i^{0}$ to be the union of $P_i^{0},\ldots,P_i^{L_i'}$ instead while still denoting it $P_i^{0}$ and renaming the $P_i^{L_i'+1},\ldots,P_i^{L_i}$ as $P_i^{1},\ldots,P_i^{L_i}$ we arrive at $|P_i^{0}|\leq\alpha|X_i|$ and $|\partial P_i^{\ell}|<\alpha|P_i^{\ell}|$ for every $1\leq\ell\leq L_i$.
\end{proof}


In addition to Corollary $\ref{almostexpanderdecomposition}$ we need a lemma that allows to transform the almost expanders into actual ones. Again the construction is due to Kun, but tweaked to work out in the bounded degree case.

\begin{Lemma}[{compare with \cite[Lemma 14]{kun}}] \label{rewire}
    Let $G$ be a finite graph with degree bounded by $d$ and $P\subseteq G$ be a subset. Assume there is $C>0$ such that for any $T\subseteq P$ with $|T|\leq |P|/2$ we have $|\partial T|\geq C|T|$. Then there is $\alpha>0$ such that if $|\partial P|<\alpha |P|$ then after adding and removing at most $\alpha |P|$ edges in $G$ and then removing at most $\frac{\alpha}{C} |P|$ vertices, we get $\partial P = \varnothing$ and $P$ becomes a connected graph with Cheeger constant at least $C/6$ and degree at most $d$.
\end{Lemma}

\begin{proof}
    The idea behind the proof is as follows: small subsets of $P$ already have expansion at least $C$ only that part of their boundary may leave $P$. If now $|\partial P|$ is small relative to $|P|$ we can amend this by rewiring these edges to sparsely chosen point inside of $P$ so that expansion is retained. For this choose $r=\frac{4}{C}$ and denote by $Q$ the set of vertices of $P$ that are endpoints of an edge in $\partial P$. If $\alpha>0$ is small enough we find a subset $F\subseteq E(P,P)$ with $|F|= |\partial P|$ such that no endpoint of an edge in $F$ is adjacent to an element of $Q$ and different edges in $F$ are at distance at least $2r$. This can be done inductively since $|B_r(A)|\leq d^{r+1}|A|$ for any $A\subseteq G$, so that as long as $|A|\leq |\partial P|$ and $\alpha<1/d^{r+1}$ we have
    \[
    |B_r(A)| < d^{r+1}\alpha |P|<|P|.
    \]

    For each $e\in F$ denote by $e^+\in P$ the endpoint of $e$ that belongs to its larger connected component obtained by removing $e$ and by $e^-$ the other endpoint of $e$. If both endpoints belong to the same connected component, we can choose either. Now take a bijection $f:\partial P\to F$ and define a new set of edges $E'(G)$ by
    \begin{align*}
        E'(G) &= \Big(E(G)\setminus \left(\partial P\cup F\right)\Big) \cup \Big\{(x,f((x,y))^+)\,:\,x\in Q, (x,y)\in \partial P\Big\}.
    \end{align*}
    If $e\in F$ is an edge such that the connected component of $e^+$ is now disjoint from that of $e^-$ then by assumption on $P$ the connected component of $e^-$ has size at most $\frac{1}{C}$. Therefore the union of all such ``new'' connected components has at most $\frac{\alpha}{C}|P|$ vertices and we now remove all of them from $P$.

    Let us now show the claim about the Cheeger constant, that is for each $T\subseteq P$ with $|T|\leq |P|/2$ we must have $|E'(T,P\setminus T)|\geq \frac{C}{6}|T|$. Let $M= \{e^+\,:\, e\in F\}$. By looking at what edges we removed from the original graph we easily obtain the inequality
    \[
    |\partial T|\leq |E'(T,P\setminus T)\cap E(T,P\setminus T)|+|T\cap M|+\sum _{a\in T\cap Q}|B_1(a)\setminus P|.
    \]
    The last term counts the edges leaving $P$ through $T$. The newly wired edges then join $T\cap Q$ either with $P\setminus T$ or $T\cap M$, so we have 
    \[
    \sum _{a\in T\cap Q}|B_1(a)\setminus P|\leq |E'(T,P\setminus T)\setminus E(T,P\setminus T)|+|T\cap M|.
    \]
    These two inequalities imply
    \[
    |\partial T|\leq |E'(T,P\setminus T)|+2|T\cap M|.
    \]
    For $x\in T\cap M$ consider the set $T\cap B_r(x)$. If $|T\cap B_r(x)|\geq r$, then $T$ contains at least $r$ elements. If however $|T\cap B_r(x)|<r$, then an endpoint of an edge in $E'(T,P\setminus T)$ must be contained in $T\cap B_r(x)$. Thus we have
    \[
    |T\cap M|\leq \frac{|T|}{r}+|E'(T,P\setminus T)|
    \]
    and with the above this gives us
    \[
    |\partial T|\leq 3|E'(T,P\setminus T)|+\frac{2|T|}{r}.
    \]
    Recall that by the choice of $T$ we have $|\partial T|\geq C|T|$ and hence we obtain
    \[
    3|E'(T,P\setminus T)|\geq |\partial T|-\frac{2|T|}{r}\geq \left(C-\frac{2}{r}\right)|T|=\frac{C}{2}|T|.
    \]
\end{proof}

Now we have all the tools to prove the main theorem of this section.

\begin{proof}[Proof of Theorem \ref{theorem: spectral gap implies almost exp}]
    Take any decreasing sequence of positive real numbers $\alpha_k\to 0$ and from Corollary \ref{almostexpanderdecomposition} find $C>0$ and $I(\alpha_k)\in \mk u$ and inductively declare $I=I_1=I(\alpha_1)$ and $I_{k+1}=I_{k}\cap I(\alpha_{k+1})$, so that $I_k\in \mk u$ since $\mk u$ is a filter. Since we have a decomposition $X_i=\bigsqcup_{\ell=0}^{L_i}P_i^\ell$ with $|P_i^{0}|<\alpha_i|X_i|$ and $|\partial P_i^\ell|<\alpha_i |P_i^\ell|$ for each $i\in I$, we can apply Lemma \ref{rewire} to every $P_i^\ell$ and throw away the $P_i^{0}$ to obtain graphs $X_i'$ which are disjoint unions of connected graph with Cheeger constant at least $C/6$. Thus $X'=\bigsqcup_iX_i'$ is an expander along $\mk u$ that is approximately isomorphic to $X$.
\end{proof}

\begin{Rem}\label{connectedinmaintheorem}
    If in the above theorem $X$ consists of connected graphs, then for every $K > 0$ the density of connected components of size at most $K$ in $X' = \bigsqcup_i X'_{i}$ converges to zero as $i\to \mk u$: if $X'_{i,(K)}$ is the union of these components, one would need at least $|X'_{i,(K)}|/K$ edges to make a connected graph out of them. Therefore, for every fixed $K>0$ we can discard all components of size at most $K$ from $X'$ retaining the approximate isomorphism to $X$. By a diagonal argument it then follows that $X'$ can be chosen such that the minimal size of a connected component in $X'_i$ converges to infinity as $i\to\mk u$. In the case that $X$ has almost boundary (T), the proof of Theorem \ref{exp+bdrDiffuse} then guarantees that $X'$ still has geometric property (T).
\end{Rem}

\section{Consequences for sofic approximations}

Having worked with general box spaces in the previous section, we will collect here the conclusions for sofic approximation.

\begin{Cor}\label{cor:sofic+T+exp}
Let $\Gamma$ be a finitely generated group with property (T) and $(X_i)$ a sofic approximation of $\Gamma$ such that each $X_i$ is connected and $X=\bigsqcup_i X_i$ is an expander. Then $X$ has geometric property (T).
\end{Cor}

\begin{proof}
    This is Proposition \ref{t+exp}.
\end{proof}

\begin{Thm}\label{thm: every sofic approx is approx T}
Let $\Gamma$ be a finitely generated group with property (T). Then every sofic approximation $(X_i)$ of $\Gamma$ is approximately isomorphic to a sofic approximation $(X'_i)$ of $\Gamma$ such that  $X'=\bigsqcup_i X'_i$ has geometric property (T).
Moreover, if each $X_i$ is connected, then $X_i'$ can also chosen to be connected.
\end{Thm}

\begin{proof}
    This follows from Theorem \ref{theorem: almost T implies appr eq to geometric T} and Remark \ref{connectedinmaintheorem}.
\end{proof}

\begin{Thm}\label{thm: almost T sofic approx}
Let $\Gamma$ be a finitely generated sofic group. The following statements are equivalent:
\begin{enumerate}
    \item $\Gamma$ has property (T),
    \item every sofic approximation of $\Gamma$ has almost boundary (T),
    \item there is a sofic approximation of $\Gamma$ with geometric property (T).
\end{enumerate}
Moreover, every sofic approximation of $\Gamma$ is approximately isomorphic to one having geometric property (T).
\end{Thm}

\begin{proof}
    The directions (i) $\Rightarrow$ (ii) and (iii) $\Rightarrow$ (i) are proved in \cite[Sec. 4]{AFS}, and (ii) $\Rightarrow$ (iii) follows from the previous theorem.
\end{proof}

Let us complete this section by considering an example of a specific sofic approximation of a property (T) group which exhibits geometric property (T) in a non-obvious manner.

\begin{Ex}
Take an expander sofic approximation $X_i'$ of $\Gamma = \SL(n,\mb Z)$, $n\geqslant 3$ (for example by quotients). Let $\Lambda$ be a finitely generated non-(T)-group with property $(\tau)$ with respect to some family of finite quotients (for example $\Lambda = \mf F_d$); take an expander box space $(Y_i)$ of $\Lambda$ such that $|Y_i|/|X_i'| \to 0$ and choose a sequence $r_i\to +\infty$ and a sequence of points $y_i\in Y_i$ 
such that $B_{r_i}(y_i)\cong B_{r_i}(1)\subset \Cay(\Lambda)$ and $|B_{r_i}(y_i)| \leq \alpha |Y_i|$ for some $0< \alpha < \frac{1}{2}$. We set $T_i\coloneqq B_{r_i}(y_i)\subset Y_i$. Set $X_i\coloneqq X'_i\sqcup Y_i$.

Take a bijection $b_i\colon Y_i\setminus T_i\to D_i\subset X_i'$ and connect $y$ to $b_i(y)$ for all $y\in Y_i\setminus T_i$. Attach a loop to all points in $X_i$ unaffected by this. This makes $X_i$ a regular graph which has degree one bigger as $X'_i$ resp. $Y_i$. Labelling the edges of $Y_i$ arbitrarily by generators of $\Gamma$ and the added edges by an additional generator which equals $1$ in $\Gamma$ obviously yields a sofic approximation of $\Gamma$ since $X_i'$ was one, so it remains to show that $X_i$ is an expander sequence.

To this end, let $K = 2C$ be the minimum of the Cheeger constants of $(Y_i)$ and $(X_i')$ and $2$. Take a subset $A_i\cup B_i\subset X_i$ such that $A_i\subset X_i'$, $B_i\subset Y_i$ and $|A_i\cup B_i| \leq |X_i|/2$. Write $\partial_{X_i'}A_i$ for the boundary of $A_i$ considered as a subset of $X_i'$ and likewise $\partial_{Y_i}B_i$. Since we have $|Y_i|/|X_i'| \to 0$, for $i$ large enough we obtain that $|\partial_{X_i'}A_i| \geq C|A_i|$. Now from our construction of the graph $X_i'$ it follows that
\[
|\partial (A_i\cup B_i)| \geq |\partial_{X_i'}A_i| + |\partial_{Y_i}B_i| + |b_i(B_i\setminus T_i)\setminus A_i|.
\]
First consider the case where $|B_i|\leq |Y_i|/2$. Then we simply have 
\[
|\partial (A_i\cup B_i)| \geq |\partial_{X_i'}A_i| + |\partial_{Y_i}B_i| \geq C(|A_i|+|B_i|).
\]
If on the other hand $|B_i|\geq |Y_i|/2$ we estimate
\[
|B_i\setminus T_i|=|B_i|-|B_i\cap T_i|\geq |B_i|-\alpha|Y_i|\geq |B_i|-2\alpha|B_i|
\]
and therefore
\brqn
    |\partial (A_i\cup B_i)| \geq |\partial_{X_i'}A_i| + |b_i(B_i\setminus T_i)\setminus A_i|\\
    \geq \frac{C}{2}|A_i| + \frac{C}{2}|b_i(B_i\setminus T_i)\cap A_i| + |b_i(B_i\setminus T_i)\setminus A_i|\\
    \geq \frac{C}{2}(|A_i|+|(B_i\setminus T_i)|)\\
    \geq \frac{C}{2}(|A_i|+|B_i|(1-2\alpha))\\
    \geq \frac{C(1-2\alpha)}{2}(|A_i|+|B_i|).
\erqn
Thus $(X_i)$ is sequence of expanders. By Proposition \ref{t+exp} it follows that $X$ has geometric property (T).

However, there is a closed invariant subset of the boundary of $X$ ``coming from the approximation of $\Lambda$''. More precisely, consider $T = \bigsqcup_i T_i$ and let $Z_T$ be the core of it; it is non-empty since $r_i\to\infty$, see Remark \ref{coreproperties}. We denote $G_T\coloneqq G(X)|_{Z_T}$. The quotient map $q_{Z_T}\colon G(X) \to G_T$ induces a homomorphism of Roe algebras
\[
q_{Z_T,X}\colon C^*_{\max} (G(X)) \to C^*_\max (G_T)
\]
On the other hand, one can consider $T$ as a subset of $Y$, the sofic approximation of $\Lambda$, hence $Z_T$ can also be considered as a closed $G(Y)$-invariant subset of $\partial\beta Y\subset \partial\beta X$. This yields another quotient map
\[
q_{Z_T,Y}\colon C^*_{\max} (G(Y)) \to C^*_\max (G_T)
\]
Since $Z_T$ is the core of subsets which are isomorphis to balls in the Cayley graph of $\Lambda$, we deduce that $G_T \cong Z_T\rtimes \Lambda$. We also observe that in this last step $Y$ could be an \emph{arbitrary} sofic approximation of $\Lambda$; this would yield the same topological dynamical system $\Lambda \lacts Z_T$. 

Moreover, one has $q_{Z_T,X}(\Delta_X) = q_{Z_T,Y}(\Delta_Y)$; in particular, this operator has spectral gap (it inherits it from the geometric property (T) of the ``big'' space $X$). We thus deduce that $Z_T$ cannot carry an invariant measure: if this were the case, the C*-algebra generated by $\Lambda$ in $C^*_\max(G_T)$ would be isomorphic on $C^*_\max(\Lambda)$ which is impossible since the Laplacian of $\Lambda$ has no spectral gap there by lack of property (T) of $\Lambda$.
\end{Ex}

It is thus of interest to understand the C*-algebra generated by $\Lambda$ in in $C^*_\max(G_T)$ since, as remarked above, it emerges from the above construction in every sofic approximation of $\Lambda$.

\begin{Prop}
Let $\Lambda$ be an exact group. Choose a sequence $r_i\to +\infty$, set $T_i\coloneqq B_{r_i}(x_i)\subset X_i\coloneqq B_{2r_i}(x_i)$ and let $Z_T$ be the core of $T \coloneqq \bigsqcup_i T_i$ in $X \coloneqq \bigsqcup_i X_i$, and set $G_T\coloneqq G(X)|_{Z_T}$. Then the C*-algebra generated by $\Lambda$ in $C^*_{\max}(G_T)\cong C(Z_T)\rtimes_\max \Lambda$ is isomorphic to $C^*_r(\Lambda)$.
\end{Prop}
\begin{proof}
Denote by $Z_\omega\coloneqq \overline{\Lambda\cdot \omega}$ the closure of the orbit of $\omega$ in $Z_T$. We claim that for every $\omega\in Z_T$ we have an isomorphism of actions $(\Lambda\lacts \beta\Lambda)\cong (\Lambda\lacts Z_\omega)$ induced by the orbit map $\lambda\mapsto \lambda\cdot\omega$. Indeed, as $Z_T$ is compact, we get an equivariant continuous map $\theta\colon \beta\Lambda\to Z_\omega$ extending the orbit map, so we only have to check is that $\theta$ is injective. Its restriction to $\Lambda$ is injective since the action on $Z_T$ is free as $r_i\to\infty$, so it remains to check injectivity on $\partial\beta\Lambda$. To this end, take two non-principal ultrafilters $\xi\neq\eta \in \partial\beta\Lambda$. Let $A,B\subset \Lambda$ be two disjoint subsets belonging to $\xi$ resp. $\eta$ and $F\in \omega$ be arbitrary; since $\omega\in Z_T$, there is a sequence $r_{f,i}\to \infty$ such that for every choice of $f_i\in F\cap T_i$ the ball $B_{r_{f,i}}(f_i)$ is isomorphic to the corresponding ball in the Cayley graph of $\Lambda$. Therefore the action of $a\in B_{r_{f,i}}(1)\subset \Lambda$ on $f_i$ is well-defined, and we set
\[
A_F\coloneqq \{a\cdot f_i\mid a\in A\in B_{r_{f,i}}(1)\cap A,\; f_i\in F\cap T_i\}\subset T
\]
and similarly for $B_F$. Since $A$ and $B$ are disjoint, so are $A_F$ and $B_F$, and therefore the clopen subsets $\overline{A_F}$ and $\overline{B_F}$ are disjoint as well. However, $\theta(\xi)\in\overline{A_F}$ and $\theta(\eta) \in \overline{B_F}$ by construction which proves that $\theta(\xi)\neq \theta(\eta)$.

Now we observe that by exactness of the maximal crossed product functor the natural map 
\[
C(Z_T)\rtimes_\max \Lambda \to \prod_{\omega\in Z_T} C(Z_\omega)\rtimes_\max\Lambda
\]
is injective. By the claim in the first paragraph we deduce that for every $\omega\in Z_T$ the crossed product $C(Z_\omega)\rtimes_\max\Lambda$ is isomorphic to $C(\beta\Lambda)\rtimes_\max\Lambda\cong \ell^\infty(\Lambda)\rtimes_\max \Lambda$. Since $\Lambda$ is assumed to be exact, the natural left action $\Lambda\lacts\ell^\infty(\Lambda)$ is amenable, so $\ell^\infty(\Lambda)\rtimes_\max \Lambda\cong \ell^\infty(\Lambda)\rtimes_r \Lambda$. In particular, the C*-algebra generated by $\Lambda$ in $\prod_{\omega\in Z_T} C(Z_\omega)\rtimes_\max\Lambda$ and hence also in $C(Z_T)\rtimes_\max \Lambda$ is isomorphic to $C^*_r(\Lambda)$, as desired.
\end{proof}

\section{A version of \.Zuk's criterion for box spaces}\label{zuksection}

In the case of property (T) for discrete groups there is a well known sufficient condition found by \.{Z}uk (see e.g. \cite{Zuk}, \cite[Chapter 5]{BHV}, \cite[Chapter 12]{brown-ozawa}): one considers the finite weighted graph with vertex set $S$ and edges as in the Cayley graph of $\Gamma$ called the \emph{link} $L(\Gamma,S)$ and only needs to inspect the smallest non-negative eigenvalue $\lambda_1$ of the graph Laplacian. The weights can be given as follows: for $s\in S$ let $\nu(s)=|\{t\in S\,:\,s\sim t\}|$ and $|\nu|=\sum_{s\in S}\nu(s)$, then the graph Laplacian becomes
\[
(\Delta^{L(\Gamma,S)} f)(s)=f(s)-\frac{1}{\nu(s)}\sum_{t\sim s}f(t)
\]
for $f\in L^2(L(\Gamma,S),\nu)$ with the scalar product given by
\[
\scp{f_1,f_2}=\sum_{s\in S}\nu(s)f_1(s)\overline{f_2(s)}.
\]

\begin{Thm}[{\.{Z}uk's criterion, see \cite[Theorem 5.6.1]{BHV}}]
    Let $\Gamma=\langle S\rangle$ be a discrete group with $S$ a finite symmetric set not containing the identity of $\Gamma$. Suppose that the link $L(\Gamma,S)$ is connected with $\lambda_1>\frac{1}{2}$. Then $\Gamma$ has property (T). Moreover, we have that $(S,\sqrt{2(2-1/\lambda_1)})$ is a Kazhdan pair for $\Gamma$.
\end{Thm}

The aim for this section is to provide a similar condition for box spaces that ensures geometric property (T), or at least almost boundary (T). We will give two ways for proving this, inspired by \cite[Theorem 12.1.15]{brown-ozawa} and \cite{ozawa-T} respectively. In both cases we use methods from \cite[Chapter 5]{BHV} where a more general version of \.Zuk's criterion is shown, that is \cite[Theorem 5.5.2]{BHV}. 

\begin{Def}
For each $x\in X$ define its \emph{linking graph} $L_x$ by
\[
L_x=B_1(x)\setminus\{x\}
\]
with edges induced from $X$, that is $E(L_x)=E(X)\cap (L_x\times L_x)$. As with the group case we need to assume that $L_x$ is connected for every $x\in X$.
\end{Def}

 Observe the following relation: we have that $(y,z)$ is an edge in $L_x$ if and only if the triangle $(x,y,z)$ lies in $X$, meaning all $(x,y)$, $(y,z)$ and $(x,z)$ are edges in $X$. The assumption that each $L_x$ be connected then implies that every edge in $X$ lies in some triangle. We will denote by $\tau(x,y)$ the number of triangles in $X$ which contain the edge $(x,y)$, and we also set $\tau(x,y)=0$ if $x$ and $y$ are identical or not adjacent at all. It follows that $\tau(x,y)=\deg_{L_x}(y)$ for the number of adjacent vertices of $y$ in $L_x$. We also set 
\[
\tau(x)=\sum_{y\sim x}\tau(x,y)
\]
which is twice the number of triangles in $X$ that contain $x$, since each such triangle gets counted for both vertices that are not $x$. 

Using this notation we will define a weighted graph structure on $L_x$. If we adopt the language of \cite[Chapter 5]{BHV} we define a ``random walk" on $L_x$ by 
\[
\mu_x(y,z)=\frac{1}{\tau(x,y)}
\]
with ``stationary measure"
\[
\nu_x(y)=\frac{\tau(x,y)}{\tau(x)}.
\]
What this means to us amounts to the following: denote by $\ell^2(L_x,\nu_x,\mathcal H)$ the Hilbert space of functions $L_x\to \mathcal H$ with the weighted scalar product 
\[
\scp{f_1,f_2}_x=\sum_{y\in L_x}\nu_x(y)\scp{f_1(y),f_2(y)}
\]
and the associated graph Laplacian 
\[
(\Delta^{L_x}f)(y)=f(y)-\sum_{y\sim z}\mu_x(y,z)f(z)
\]
which is then a positive operator on $\ell^2(L_x,\nu_x,\mathcal H)$. Writing $\lambda_1(L_x)$ for the first positive eigenvalue of $\Delta^{L_x}$ we have the Poincar\'e inequality
\[
\lambda_1(L_x) \sum_{y,z\in L_x}\norm{f(y)-f(z)}^2\nu_x(y)\nu_x(z) \leq \sum_{y,z\in L_x}\norm{f(y)-f(z)}^2 \nu_x(y)\mu_x(y,z),
\]
see \cite[Proposition 5.3.1]{BHV}.

In order to obtain a connection between these weighted graphs and geometric property (T) we need to introduce some related operators in $\mathbb C[X]$. Define $\Delta_{\tau}\in\mathbb C[X]$ by 
\[
(\Delta_{\tau}\eta)(x)=\sum_{y\sim x}(\eta(x)-\eta(y))\tau(x,y)
\]
for $\eta\in \ell^2 X$. If we also set 
\[
(N_{\tau}\eta)(x)=\tau(x)\eta(x),\quad (M_{\tau}\eta)(x)=\sum_{y\sim x}\tau(x,y)\eta(y)
\]
then $\Delta_{\tau}=N_{\tau}-M_{\tau}$. Our aim will be to proof spectral gap for $\Delta_\tau$ and this will be enough to conclude geometric property (T), as the following lemma will show.

\begin{Lemma}
    Let $X$ be a box space of graphs such that for each $x\in X$ the linking graph $L_x$ is connected. Then $X$ has geometric property (T) if and only if $\Delta_\tau$ has spectral gap in $C^*_{u,\max}(X)$. More precisely, spectral gap for $\Delta$ is equivalent to spectral gap for $\Delta_\tau$.
\end{Lemma}

\begin{proof}
    As was mentioned in a previous section, it is proven in \cite[Section 5]{WYgeomT} that $\Delta$ can be written as a finite sum 
    \[
    \Delta=\sum_{j=1}^n\Delta^{E_j}
    \]
    where $\Delta^{E_j}=(v_jv_j^*-v_j)^*(v_jv_j^*-v_j)$ and the $v_j\in \mathbb C[X]$ are partial isometries whose source and range are disjoint and the $E_j$ are the graphs of the partial translation underlying the $v_j$. By decomposing these partial isometries further we can assume that $\tau(\cdot,\cdot)$ takes constant values $\alpha_j$ along the graphs of each partial translation. This is possible because $\tau(\cdot,\cdot)$ takes only finitely many values. Moreover, as discussed in the above paragraph, every linking graph being connected implies $1\leq \tau(\cdot,\cdot)\leq d$, where $d$ is again the uniform bound on the degree of $X$. We conclude
    \[
    \Delta_\tau = \sum_{j=1}^n \alpha_j\Delta^{E_j}
    \]
    and observe that for any representation $\pi:\mathbb C[X]\to\mathcal B(\mathcal H)$ and $\xi\in \mathcal H$ we have
    \[
    \scp{\pi(\Delta)\xi,\xi}\leq\scp{\pi(\Delta_\tau)\xi,\xi}\leq d\scp{\pi(\Delta)\xi,\xi}.
    \]
    From this the assertion follows easily.
\end{proof}

\begin{Thm}\label{thm:geomZuk}
    Let $X=\bigsqcup_i X_i$ be a box space of graphs such that each linking graph $L_x$ is connected. Suppose that there are $Y_i\subseteq X_i$ with $|Y_i|/|X_i|\to 1$ and $\lambda > \frac{1}{2}$ such that for every $x\in Y_i$ we have 
    \[
    \lambda_1\left(L_x\right)\geq \lambda.
    \]
    Then $X$ has almost boundary (T) with respect to $\mathrm{core}(Y)$. If $Y_i = X_i$ for all $i\in\mathbb N$, then $X$ has geometric property (T).
\end{Thm}

To show spectral gap for $\Delta_{\tau}$ in this setting we will directly write the operator $\Delta_{\tau}^2-c\Delta_{\tau}$ as a sum of squares in $\mathbb C[X]$, that is
\[
\Delta_{\tau}^2-c\Delta_{\tau}=\sum_j S_j^*S_j
\]
for some $S_j\in\mathbb C[X]$. That way, spectral gap of $\Delta_\tau$ in $C_{u,\max}^*(X)$ is immediate. We will achieve this by employing ideas of \cite{ozawa-T} concerning property (T) for groups, but again, applied to to $L_x$ instead of $L(\Gamma,S)$. Specifically, the Poincar\'{e} inequality for $L_x$ yields (see \cite[Example 5]{ozawa-T}) that $\lambda^{-1}\Delta^{L_x} + P_x - I_x $ is a positive operator on $\ell^2(L_x,\nu_x)$, where $P_x$ denotes the projection onto the subspace of constant functions. Hence we have
\[
\lambda^{-1}\Delta^{L_x} + P_x - I_x = T_x^*T_x
\]
for some $T_x\in \mathcal B(\ell^2(L_x,\nu_x))$. 
\begin{proof}[Proof of Theorem \ref{thm:geomZuk}]
    Note that each $B_1(x)$ is a graph on at maximum $d+1$ vertices, where $d$ is the uniform bound on the degree of $X$. Hence there are only finitely many graphs $G_1,\ldots,G_n$ that appear as some $B_1(x)$ and consider them as rooted graphs, that is with the vertex corresponding to $x$ marked as a root. Take now a partition of $X$ into disjoint subsets $C_1,\ldots,C_n$ such that we have graph isomorphisms
    \[
    \varphi_{j,x}: G_j\to B_1(x),\quad x\in C_j
    \]
    for $j=1,\ldots,n$. We may assume that for each $j$ the images $(\varphi_{j,x}(G_j))_{x\in C_j}$ are disjoint in $X$ if we allow repetitions among the $G_1,\ldots,G_n$. Similarly we can and will assume that the union of the first $m\leq n$ subsets $C_1,\dots,C_m$ equals $Y = \bigcup_i Y_i$. By removing the root from each $G_j$ we obtain graphs $L_1,\ldots,L_n$ with graph isomorphisms
    \[
    \varphi_{j,x}: L_j\to L_x,\quad x\in C_j
    \]
    for $j=1,\ldots,n$. As in the above paragraph, we want to equip our graphs $L_j$ with the weights $\nu_j(\ell)=\frac{\deg_{L_j}(\ell)}{\tau_j}$ and $\mu_j(\ell,\ell')=\frac{1}{\deg_{L_j}(\ell)}$ where $\tau_j=\sum_{\ell\in L_j}\deg_{L_j}(\ell)$. Let us further define for every $\ell\in L_j$ a partial bijection
    \begin{align*}
        t_{j,\ell}: C_j &\to t_{j,\ell}(C_j)\\
        x &\mapsto \varphi_{j,x}(\ell)
    \end{align*}
    and denote by $v_{j,\ell}\in\mathbb C[X]$ the corresponding partial isometry. All of this is to obtain a description of the Laplace operator on $X$ that is sensible to the information of each $L_x$. 

    We adopt the notation from the first paragraph of the first proof. The Poincar\'{e} inequality for $L_j$ is then equivalent to the following identity in $\mathcal B(\ell^2(L_j,\nu_j))$:
    \[
    \lambda^{-1} \Delta^{L_j} + P_j - I_j = T^*_j T_j   
    \]
    for some $T_j\in\mathcal B(\ell^2(L_j,\nu_j))$, and we need to convert this 'local' information into 'global' information expressed in terms of our previously defined operators $\Delta_{\tau}$, $N_{\tau}$ and $M_{\tau}$. Similar to \cite[Example 5]{ozawa-T} consider the element $v_j$ given by
    \[
    v^*_j = \sum_{\ell\in L_j} \delta_\ell\otimes  v^*_{j,\ell} \in \ell^2(L_j,\nu_j)\otimes \mathbb C[X],
    \]
    Note that for every operator $T\in \mathcal B(\ell^2(L_j,\nu_j))$ we have a $\mathbb C[X]$-linear endomorphism $T\otimes I$ of $\ell^2(L_j,\nu_j)\otimes \mathbb C[X]$, so that we can compute its $\mathbb C[X]$-valued matrix coefficient 
    \[
    \scp{v_j^*,T v_j^*}_{\mathbb C[X]}=\sum_{\ell,\ell'\in L_j}\scp{\delta_\ell,T\delta_{\ell'}}_j v_{j,\ell}v_{j,\ell'}^*.
    \]
    Thus, applying $\scp{v_j^*,\,\cdot\, v_j^*}_{\mathbb C[X]}$ to both sides of
    \[
    \lambda^{-1} \Delta^{L_j} + P_j - 1_j = T^*_j T_j,
    \]
    we obtain on the right-hand side $\scp{T_j v_j^*, T_j v_j^*}_{\mathbb C[X]}$, which is by construction a sum of squares in $\mathbb C[X]$. On the left-hand side we get the sum of three terms, which calculate as follows:
    \begin{align*}
        \lambda^{-1} \scp{v^*_j,\Delta^{L_j} v_j^*}_{\mathbb C[X]} &= \lambda^{-1}\sum_{\ell\in L_j} \nu_j(\ell)v_{j,\ell}v_{j,\ell}^* - \lambda^{-1}\sum_{\ell\sim\ell'\in L_j} \frac{1}{\tau_j} v_{j,\ell}v_{j,\ell'}^*,\\
        \scp{v^*_j,P_jv_j^*}_{\mathbb C[X]} &= \sum_{\ell,\ell'\in L_j} \nu_j(\ell)\nu_j(\ell') v_{j,\ell}v_{j,\ell}^*,\\
        -\scp{v^*_j,v_j^*}_{\mathbb C[X]} &= -\sum_{\ell\in L_j} \nu_j(\ell) v_{j,\ell} v_{j,\ell}^*.
    \end{align*}
    Multiplying by $\tau_j$ and summing these equations over $j = 1,\ldots,n$ (or $1,\ldots,m$ respectively) we of course still get a sum of squares on the right-hand side. On the left-hand side we get another sum of three terms. The first one is
    \begin{align*}
        \sum_{j=1}^n \lambda^{-1} \tau_j\scp{v^*_j,\Delta^{L_j} v_j^*}_{\mathbb C[X]}
        &= \lambda^{-1}\sum_{j=1}^n\sum_{\ell\in L_j} \deg_{L_j}(\ell)v_{j,\ell}v_{j,\ell}^* - \lambda^{-1}\sum_{j=1}^n\sum_{\ell\sim\ell'\in L_j} v_{j,\ell}v_{j,\ell'}^*\\
        &=\vphantom{\sum_{j=1}^n} \lambda^{-1}(N_{\tau}-M_{\tau}) = \lambda^{-1}\Delta_{\tau}.
    \end{align*}
    The second term is
    \begin{align*}
    \sum_{j=1}^n \tau_j \scp{v^*_j,P_jv_j^*}_{\mathbb C[X]} &= \sum_{j=1}^n \sum_{\ell,\ell'\in L_j} \tau_j\nu_j(\ell)\nu_j(\ell') v_{j,\ell}v_{j,\ell}^*\\
    &= \sum_{k=1}^n\sum_{p=1}^n \sum_{m\sim m'\in L_k}\sum_{q\sim q'\in L_p} v_{k,m}v^*_{k,m'} N_\tau^{-1} v_{p,q}v^*_{p,q'}\\
    &= \vphantom{\sum_{j=1}^n} M_{\tau} N_{\tau}^{-1} M_{\tau}.
    \end{align*}
    To illustrated the second to last equality consider the following diagram:
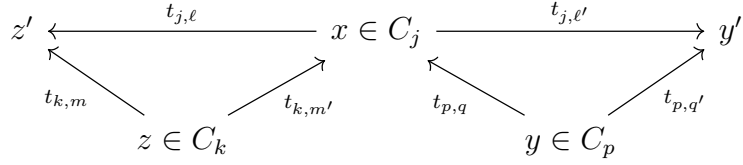
\begin{center}
\begin{tikzcd}
z' &                                                            & x\in C_j \arrow[ll, "{t_{j,\ell}}"'] \arrow[rr, "{t_{j,\ell'}}"] &                                                            & y' \\
   & z\in C_k \arrow[ru, "{t_{k,m'}}"'] \arrow[lu, "{t_{k,m}}"] &                                                                 & y\in C_p \arrow[ru, "{t_{p,q'}}"'] \arrow[lu, "{t_{p,q}}"] &   
\end{tikzcd}
\end{center}
    Note that every arrow $t_{j,\ell}$ is counted $\tau_j\nu_j(\ell)$ times, which is the number of triangles containing that arrow. The same holds for $t_{j,\ell'}$ with $\tau_j\nu_j(\ell')$, and since applying $N_{\tau}^{-1}$ is just multiplying with $\tau_j^{-1}$ at $x\in C_j$, one of the factors $\tau_j$ vanishes. Lastly, the third term computes as
    \begin{align*}
    -\sum_{j=1}^n \tau_j \scp{v^*_j,v_j^*}_{\mathbb C[X]}=-\sum_{j=1}^n \sum_{\ell\in L_j} \tau_j \nu_j(\ell) v_{j,\ell} v_{j,\ell}^* = -N_\tau.
    \end{align*}
    
    Putting these together we see that
    \[
    \lambda^{-1}\Delta_{\tau} + M_\tau N_\tau^{-1} M_\tau - N_\tau
    \]
    is a sum of squares in $\mathbb C[X]$. As in the first proof, this is equal to
    \[
    \Delta_\tau - \lambda(2\Delta_\tau-\Delta_\tau N_{\tau}^{-1}\Delta_\tau)
    \]
    and since $0\leq I-N_{\tau}^{-1}\leq 1$ is itself a square in $\mathbb C[X]$, we get that for $c = 2 - 1/\lambda$
    \[
    \Delta_{\tau}^2 - c\Delta_\tau
    \]
    is a sum of squares in $\mathbb C[X]$. As explained before, this means spectral gap at least $c>0$ in $C_{u,\max}^*(X)$ (or $C^*(\partial G(X)|_{\mathrm{core}(Y)})$ respectively).
\end{proof}

\begin{Rem}
    The spectral gap $c=2-1/\lambda$ in the above proof is the same term appearing in the Kazhdan constant from \.Zuk's theorem $\varepsilon=\sqrt{2c}$. This is of course no coincidence, but a consequence of the Poincar\'e inequality coming from the local information of the graphs.
\end{Rem}



\section{Concluding remarks and open questions}

We finish the paper by some comments and questions which seem natural in view of the above results.

First of all, our main result gives many examples of box spaces with almost boundary (T) coming from Farber sequences in property (T) groups. Recall that a sequence of finite index subgroups $\Gamma_i$ in $\Gamma$ is called Farber if the box space of Schreier graphs $X_i = \Gamma/\Gamma_i$ is a sofic approximation of $\Gamma$. Since it is known that all sequences of finite index subgroups $\Gamma_i$ in a higher rank lattice $\Gamma$ with $[\Gamma:\Gamma_i]\to \infty$ are Farber \cite{7samurai,abert-g-n}. Since $\Gamma$ acts on $X_i = \Gamma/\Gamma_i$, the corresponding box space is a (connected) expander, and by Corollary \ref{cor:sofic+T+exp} it has geometric property (T). However, geometric property (T) in this case remains inherited from property (T) of the group, so it would be interesting to get a ``group-independent'' construction.

\begin{Qu}
Is there a randomized construction of geometric property (T) box spaces of graphs in the spirit of \,{Z}uk's construction for groups?
\end{Qu}

Since Theorem \ref{thm:geomZuk} is based on representing the Laplace operator as a sum of squares, it is natural to ask whether it characterizes geometric property (T):

\begin{Qu}
Is geometric property (T) equivalent to $\Delta^2-c\Delta$ being a sum of squares in $\mb C[X]$?
\end{Qu}

The comparison of topological and measured versions of property (T) also leave some open questions about their connections.

\begin{Qu}
Is there a box space of graphs with boundary geometric property (T) which is not an expander (equivalently, does not have geometric property (T))?
\end{Qu}

\begin{Qu}
Let $X$ be a box space of graphs such that its coarse boundary groupoid has has measured property (T) for every invariant measure on $\partial\beta X$. Does it necessarily have geometric property (T)?
\end{Qu}
The corresponding natural question for étale topological groupoids would be whether having measured property (T) for every quasi-invariant measure on the base space implies topological property (T) in general; however, this would require to study the relation between these notions beyond the p.m.p. case which we focused on in this paper.

Another interesting property of a box space is the type of the measure algebra of invariant sets in $(\partial\beta X,\mu_{\mk u})$; let us call a box space \emph{diffuse} resp. \emph{atomic} if this measure algebra is of the corresponding type. Our main result characterizing spectral gap (Corollary \ref{corollary:almost T is appr expander}) is particularly interesting when the expander decomposition is diffuse; in fact, in the opposite atomic case, one can obtain a much simpler proof since every atom necessarily gives an expander. The discrepancy between the diffuse and the atomic case is at heart of the results of Kun and Thom providing first examples of non-amenable groups which don't admit a sofic approximation by a (connected) expander \cite{kun-thom, kun-youtube} 
The following is a generalized version of one of the crucial steps in their construction:

\begin{Prop}[{see \cite{kun-youtube}}]
Let $N\trianglelefteq\Gamma$ be a normal subgroup such that the quotient $Q = \Gamma/N$ is amenable. If $\Gamma$ has a sofic approximation by a connected expander, then its restriction to $N$ is atomic; in particular, $N$ itself admits an approximation by a connected expander.
\end{Prop}
\begin{proof}
Let $X$ be the corresponding box space. Since it is an expander, $\Gamma$ has spectral gap in the ultraproduct Hilbert space $\leftsub{\mk u}{\mc H}$ and its action on $(\partial\beta X,\mu_{\mk u})$ is ergodic; these properties are then retained by the action of $Q$ on the subspace $\leftsub{\mk u}{\mc H}^N$ of $N$-invariant vectors. But since $Q$ is amenable, an ergodic action with spectral gap is discrete (which by ergodicity implies that $L^\infty(\partial\beta X,\mu_{\mk u})^N$ is finite-dimensional). Each atom in $L^\infty(\partial\beta X,\mu_{\mk u})^N$ gives a minimal $N$-invariant subset, and thus yields a connected expander in the restriction of the sofic approximation to $N$.
\end{proof}

In \cite{kun-youtube} this fact is then played against the main result of \cite{kun-thom} which allows to conclude that the centralizer of $N$ in $\Gamma$ must be LEF; it follows that in every sofic approximation of the group $\Gamma = N\times\Delta$ where $N$ has property (T) and $\Delta$ is amenable non-LEF, the restriction of the approximation to $N$ necessarily produces a \emph{diffuse} expander. It is therefore natural to ask the following:
\begin{Qu}
Let $\Lambda\leq \Gamma$ be an inclusion of finitely generated groups. Under which conditions does every sofic approximation of $\Gamma$ restrict to a diffuse sofic approximation of $\Lambda$?
\end{Qu}

Finally, it is natural to ask whether the result concerning the spectral gap in the ultraproduct has a natural analogue in the hyperlinear case: 
\begin{Qu}
Let $\Gamma = \langle S\rangle$ be a finitely generated hyperlinear group. Consider a hyperlinear approximation of $\Gamma$, i.e. an embedding of $\Gamma$ into the the metric ultraproduct $\prod_{n\to\mk u} U(n)$ of unitary matrices. Suppose that the Laplace operator of $\Gamma$ 
\[
\Delta_\Gamma = \sum_{s\in S} (1-s)^*(1-s)
\]
has spectral gap when acting on the ultraproduct $\leftsub{{n\to\mk u}}{\prod} \mb C^n$. Can one then realise the given hyperlinear approximation of $\Gamma$ by tuples of unitary matrices $(u_1^{(n)},\dots,u_k^{(n)})\in U(n)^k$ whose Laplacians
\[
\Delta_n = \sum_{j=1}^k \left(1-u_{j}^{(n)}\right)^*\left(1-u_{j}^{(n)}\right)
\]
have uniform spectral gap along $\mk u$?
\end{Qu}

\newcommand{\etalchar}[1]{$^{#1}$}


\begin{thebibliography}{ABB{\etalchar{+}}17}

\bibitem[AB21]{alekseev-biz}
Vadim Alekseev and Leonardo Biz.
\newblock Sofic boundaries and a-{{T-menability}}.
\newblock (arXiv:2112.14760), 2021.

\bibitem[ABB{\etalchar{+}}17]{7samurai}
Miklos Abert, Nicolas Bergeron, Ian Biringer, Tsachik Gelander, Nikolay Nikolov, Jean Raimbault, and Iddo Samet.
\newblock On the growth of {$L^2$}-invariants for sequences of lattices in {L}ie groups.
\newblock {\em Ann. of Math. (2)}, 185(3):711--790, 2017.

\bibitem[AD05]{A-D}
Claire Anantharaman-Delaroche.
\newblock Cohomology of property {$T$} groupoids and applications.
\newblock {\em Ergodic Theory Dynam. Systems}, 25(4):977--1013, 2005.

\bibitem[ADR00]{amenable-groupoids}
Claire Anantharaman-Delaroche and Jean Renault.
\newblock {\em Amenable groupoids}, volume~36 of {\em Monographies de L'Enseignement Math\'ematique}.
\newblock L'Enseignement Math\'ematique, Geneva, 2000.
\newblock With a foreword by Georges Skandalis and Appendix B by E. Germain.

\bibitem[AF19]{AFS}
Vadim Alekseev and Martin {Finn-Sell}.
\newblock Sofic boundaries of groups and coarse geometry of sofic approximations.
\newblock {\em Groups Geom. Dyn.}, 13(1):191--234, 2019.

\bibitem[AGN17]{abert-g-n}
Miklos Abert, Tsachik Gelander, and Nikolay Nikolov.
\newblock Rank, combinatorial cost, and homology torsion growth in higher rank lattices.
\newblock {\em Duke Math. J.}, 166(15):2925--2964, 2017.

\bibitem[BdV08]{BHV}
Bachir Bekka, Pierre {de la Harpe}, and Alain Valette.
\newblock {\em Kazhdan's Property ({{T}})}, volume~11 of {\em New {{Mathematical Monographs}}}.
\newblock Cambridge University Press, Cambridge, 2008.

\bibitem[BO08]{brown-ozawa}
Nathanial~P. Brown and Narutaka Ozawa.
\newblock {\em {{C}}*-Algebras and Finite-Dimensional Approximations}, volume~88 of {\em Graduate {{Studies}} in {{Mathematics}}}.
\newblock American Mathematical Society, Providence, RI, 2008.

\bibitem[Con76]{connes76}
Alain Connes.
\newblock Classification of injective factors. {C}ases {$II\sb{1},$} {$II\sb{\infty },$} {$III\sb{\lambda },$} {$\lambda \not=1$}.
\newblock {\em Ann. of Math. (2)}, 104(1):73--115, 1976.

\bibitem[DW22]{dA-W}
Cl{\'e}ment Dell'Aiera and Rufus Willett.
\newblock Topological property ({{T}}) for groupoids.
\newblock {\em Ann. Inst. Fourier (Grenoble)}, 72(3):1097--1148, 2022.

\bibitem[Kai19]{kaiser}
Tom Kaiser.
\newblock Combinatorial cost: A coarse setting.
\newblock {\em Trans. Amer. Math. Soc.}, 372(4):2855--2874, 2019.

\bibitem[KT19]{kun-thom}
Gabor Kun and Andreas Thom.
\newblock Inapproximability of actions and {{Kazhdan}}'s property ({{T}}).
\newblock (arXiv:1901.03963), 2019.

\bibitem[Kun19]{kun}
Gabor Kun.
\newblock On sofic approximations of {{Property}} ({{T}}) groups.
\newblock (arXiv:1606.04471), 2019.

\bibitem[Kun21]{kun-youtube}
Gabor Kun.
\newblock Non-amenable groups admitting no sofic approximation by expander graphs.
\newblock \url{https://www.youtube.com/watch?v=jAiV6QHfcy4}, February 2021.

\bibitem[Oza16]{ozawa-T}
Narutaka Ozawa.
\newblock Noncommutative real algebraic geometry of {K}azhdan's property ({T}).
\newblock {\em J. Inst. Math. Jussieu}, 15(1):85--90, 2016.

\bibitem[Pes08]{Pestov}
Vladimir~G. Pestov.
\newblock Hyperlinear and sofic groups: a brief guide.
\newblock {\em Bull. Symbolic Logic}, 14(4):449--480, 2008.

\bibitem[Roe03]{Roe}
John Roe.
\newblock {\em Lectures on Coarse Geometry}, volume~31 of {\em University {{Lecture Series}}}.
\newblock American Mathematical Society, Providence, RI, 2003.

\bibitem[Sau05]{sauer}
Roman Sauer.
\newblock {{L}}{$^2$}-{{Betti}} numbers of discrete measured groupoids.
\newblock {\em Int. J. Algebra Comput.}, 15(05n06):1169--1188, October 2005.

\bibitem[Win24]{Winkel}
Jeroen Winkel.
\newblock Cycles in graphs with geometric property ({T}).
\newblock {\em Groups Geom. Dyn.}, 18(1):361--377, 2024.

\bibitem[WY14]{WYgeomT}
Rufus Willett and Guoliang Yu.
\newblock Geometric property ({{T}}).
\newblock {\em Chin. Ann. Math. Ser. B}, 35(5):761--800, 2014.

\bibitem[\.Z02]{Zuk}
Andrzej \.Zuk.
\newblock On property ({T}) for discrete groups.
\newblock In {\em Rigidity in dynamics and geometry ({C}ambridge, 2000)}, pages 473--482. Springer, Berlin, 2002.

\end{thebibliography}
\end{document}